\input amstex          
 \documentstyle{amsppt}  

 
\magnification=1100
\hsize=6.5truein
\vsize=8.9truein
\nologo  
\NoBlackBoxes
\TagsOnRight
\CenteredTagsOnSplits
\leftheadtext{ }
\rightheadtext{ }

\parskip=3pt
\font\first=cmsl8

\document

\baselineskip=12pt  
\parindent=20pt
\parskip=3pt


\def\A{\Cal A}
\def\C{\Cal C}
\def\DD{\Cal D}
\def\E{\Cal E}
 
\def\O{\Cal O}
\def\T{\Cal T}

\def\ra{\longrightarrow}
\def\rA#1{\overset {#1}\to\longrightarrow}
\def\ds#1{\textstyle\bigoplus\limits_{#1}}

\def\sgm#1{\textstyle\sum\limits_{#1}}

\def\cupp{{\raise 0.9pt\hbox{$\,\scriptstyle\cup\,$}}}
\def\hmat#1#2{{\raise 0.7pt 
\hbox{$\scriptstyle(#1 \hskip .04in #2)$}}}

\def\a{\alpha}
\def\b{\beta}
\def\c{\gamma}
\def\D{\Delta}
\def\d{\delta}
\def\e{\epsilon}
\def\g{\gamma}
\def\l{\lambda}
\def\p{\varphi}
\def\r{\rho}
\def\s{\sigma}
\def\th{\theta}
\def\w{\omega}
\def\wt{\widetilde}

\def\qq{\Bbb Q}

\def\zz{\Bbb Z}

\def\sd{\sqrt{\Delta}}
\def\sc{\sqrt{c}} 
\def\sm{\sqrt{-3}}
\def\cbrt{{}^3\!\!\!\!\sqrt}
\def\({\big (} 
\def\){\big )} 

\def\dd{\!\cdot\!}
\def\oo{\!\circ\!}

\def\wt{\widetilde}
\def\ol{\overline}

\def\frp{\frak p}
\def\frq{\frak q}
\def\vp{v_{\frak p}}

\def\ok{\otimes_k}

\def\rk{g}

\def\Gal{\Cal G \text{\it al\,}}
\def\kp{\widehat{k}_{\frp}}
\def\kq{\widehat{k}_{\frq}}
\def\k#1{\widehat{k}_{#1}}
\def\q#1{\widehat{\qq}_{#1}}

\def\Br{\text{\sl Br}}      
\def\br{\text{\first Br}}      

\def\char{\text{\sl char}}

\def\cor{\text{\sl cor}}
\def\coker{\text{\sl coker}}

\def\det{\text{\sl det}}
 
\def\dim{\text{\sl dim}}
\def\disc{\text{\sl disc}}

\def\id{\text{\sl id}} 
\def\im{\text{\sl im}}
\def\fim{\text{\first im}}

\def\inv{\text{\sl inv}}

\def\Hom{\text{\sl Hom}}
\def\ker{\text{\sl ker}}
\def\fker{\text{\first ker}}
\def\lim{\text{\sl lim}}

\def\min{\text{\sl min}}
\def\mod{\text{\sl mod}}

\def\rank{\text{\sl rank}}
\def\res{\text{\sl res}}
\def\fres{\text{\first res}}

\def\tors{\text{\first tors}}

\def\fh{\widehat f}
\def\ch{\widehat c}
\def\dh{\widehat d}
\def\rh{\widehat r}
\def\sh{\widehat s}
\def\xh{\widehat x}
\def\yh{\widehat y}
\def\Eh{\widehat \E}   
\def\Ehp{\Eh^{\,'}}
\def\Th{\widehat \T}
\def\Thp{\Th^{\,'}}  
\def\Kh{\widehat K}   
\def\uh{\widehat u}   

\def\afk{\text{\it A}_{f,k}}
\def\afhk{\text{\it A}_{\fh,k}}
\def\bfk{\text{\it Q}_{f,k}}
\def\bfhk{\text{\it Q}_{\fh,k}}
\def\afl{\text{\it A}_{f,L}}
\def\afhl{\text{\it A}_{\fh,L}}
\def\bfl{\text{\it Q}_{f,L}}
\def\bfhl{\text{\it Q}_{\fh,L}}
\def\aff{\text{\it A}_{f,F}}
\def\bff{\text{\it Q}_{f,F}}
\def\z{\overline{z}}
\def\M#1#2{\!\left(\!\smallmatrix #1 \\ #2 \endsmallmatrix\!\right)}	

\def\kt{\widetilde k}
\def\Lt{\widetilde L}
\def\Ft{\widetilde F}

\def\rf{\frak r}
\def\sf{\frak s}
\def\rfh{\widehat{\frak r}}
\def\sfh{\widehat{\frak s}}
\def\rfhp{\widehat{\frak r}\hskip .01 in'}
\def\sfhp{\widehat{\frak s}\hskip .01 in'}


\topmatter

\title  Curves $\C$  that are cyclic twists of $Y^2=X^3+c$ and 
the relative Brauer groups $\Br (k(\C)/k)$                      \endtitle

\thanks{We would like to thank D.~Krashen for some useful 
conversations, and in particular for his assistance with the 
proof of Theorem  6.1.}
\endthanks

\author Darrell E. Haile, \ Ilseop Han, \ Adrian R. Wadsworth      \endauthor  

\address{\it Department of Mathematics, Indiana University, 
Bloomington, IN 47405}\endaddress

\address{\it  Department of Mathematics, 
California State University, San Bernardino, CA 92407}
\endaddress

\address{\it Department of Mathematics, University of California, San Diego, 
La Jolla, CA 92093-0112}   \endaddress

\abstract  \nofrills   Let $k$ be a field with $\eightsl{char}(k) 
\ne 2,3$.
  Let $\C_f$ be the projective curve of a binary cubic form $f$,
and $k(\C_f)$ the function field of  $\C_f$.  In this paper 
we explicitly describe the relative Brauer group $\br(k(\C_f)/k)$ of 
$k(\C_f)$ over $k$.    When $f$ is 
diagonalizable we show that every algebra in $\br(k(\C_f)/k)$ 
is a cyclic algebra obtainable using
the $y$-coordinate of a $k$-rational point on the 
Jacobian $\E$ of  $\C_f$.  
But when $f$ is not diagonalizable, the algebras in 
$\br(k(\C_f)/k)$ are presented as cup products of 
cohomology classes, but not as cyclic algebras.
In particular, we provide several 
specific examples of 
relative Brauer groups for $k=\qq$, the rationals, and for 
$k=\qq(\w)$ 
where $\w$ is a primitive third root of unity.  
The approach is to realize $\C_f$ as a cyclic twist of 
its Jacobian $\E$,  an elliptic curve, and then 
apply a recent theorem
of Ciperiani and Krashen.
 \endabstract
\endtopmatter

\document
 

\noindent {\bf \S1. \ Introduction}            \medskip


Let $k$ be a field of characteristic not $2$ or $3$, and let 
$\C_f$ be the 
smooth projective curve over $k$ given by the equation 
$Z^3 = f(X,Y)$, where 
$$
f(X,Y)  \ = \ AX^3 +3BX^2Y + 3CXY^2 +DY^3 \in k[X,Y]
\tag 1.1
$$
is a nongedenerate binary cubic form.  Let $k(\C_f)$ be the function 
field of 
$\C_f$ over $k$.  The goal of this paper is to compute the relative 
Brauer group
$\Br(k(\C_f)/k)$, i.e., the kernel of the scalar extension map 
$\Br(k) \to \Br(k(\C_f))$.  The curve $\C_f$ has genus $1$, and 
is therefore  a 
principal homogeneous space of  the  elliptic curve $\E$ which is 
the Jacobian of $\C_f$. 
It was known long ago by results of Lichtenbaum in [L$_1$, 
Lemmas 1, 2, 3], [L$_2$, \S 2] 
that since $\C_f$ satisfies period = index, the elements of 
$\Br(k(\C_f)/k)$ are parametrized by $\E(k)$, the group of 
$k$-rational points 
of~$\E$.  But Lichtenbaum's mapping is difficult to use for 
explicit calculations.  The main tool we apply here is a 
wonderful recent theorem of Ciperiani and Krashen
[CK, Th.~2.6.5], which gives a cup product formula for elements of 
$\Br(k(\C)/k)$
when $\C$ is what they call a \lq\lq cyclic twist" of an elliptic 
curve $\E$.
That is, $\C$ is 
a principal homogeneous space of  $\E$ 
such that $\C$ lies in the image of $H^1(k, \T) \to H^1(k, \E)$ for 
some Galois stable cyclic subgroup $\T$ of $\E(k_s)$.  
  We prove here that for appropriate choices of $\E$ and $\T$,
 the principal homogeneous spaces over $\E$ that arise from
$H^1(k,\T)$ are precisely the curves~$\C_f$ as above, and we use 
their result to describe $\Br(k(\C_f)/k)$.

Our interest in these particular curves $\C_f$ comes in part from 
the theory of 
generalized Clifford algebras.  For every homogeneous form 
$f(X_1, \ldots, X_n) \in 
k[X_1, \ldots, X_n]$ of degree $d$ (say),\break one can define a 
$k$-algebra $A_{f,k}$, the 
Clifford algebra of $f$ over $k$, to be the quotient\break of the
free 
associative 
algebra $k\{u_1, \ldots, u_n\}$ modulo the ideal generated by the 
set\break
$\{(a_1u_1+ a_2u_2+\ldots +a_nu_n)^d - f(a_1, \ldots, a_n) \mid
a_1, \ldots, a_n\in k\}$.   In the case where $f$ is a binary 
cubic form
as in (1.1), it was shown by Haile in [H$_{1}$, Th.1.1$'$] that 
$A_{f,k}$ is an Azumaya algebra of rank $9$ over its center, and 
that its center is 
the coordinate ring of the affine curve $Y^2= X^3 - 27 \Delta_f$ 
where $\Delta_f$ is the discriminant of~$f$ (see (6.3) below 
for a formula for $\Delta_f$).  It follows that the $k$-rational 
points on this curve give rise to simple images of $A_{f,k}$ with 
center $k$.  In the same paper Haile also showed that if
we let $\E_f$ denote the projective closure of 
this affine curve, then 
the resulting function from the group of $k$-rational points on 
$\E_f$ to $\Br(k)$ (sending the point at infinity to $0$) is a 
group homomorphism,  and the image of the homomorphism is  
the relative Brauer group $\Br(k(\C_f)/k)$, where 
$\C_f $~is the curve given above.  This point of view can be used 
to compute $\Br(k(\C_f)/k)$ in many cases, such as when $k$ is 
perfect and $\Delta, -3\in k^{*2}$. (See the remarks in the proof 
of Proposition 4.1.)  The algebras in the relative Brauer group are all 
symbol algebras of degree $3$  with one component fixed, and the 
other depending on points on $\E_f(k)$.
However, the machinery of Ciperiani and Krashen has wider 
applicability.
Moreover, it is worthwhile to understand the connection between 
the two approaches.

Here is a description of the theorem of Ciperiani and Krashen that 
we will be 
using below:  Let~$\E$~be a smooth projective elliptic curve over $k$, 
and let 
$\T$ be a finite cyclic subgroup of the group~$\E(k_s)$ of points
of $\E$ over the separable closure $k_s$ of $k$.  Suppose 
$\char(k) \nmid |\T|$ and $\T$ is setwise 
stable under the action of the absolute Galois group $G_k
= \Gal (k_s/k)$.  Let $\E'$ be the elliptic curve isogenous to $\E$
such that $\T$ is the kernel of an isogeny $\lambda\: \E\to \E'$.
Let $\T'$ be the kernel of the dual isogeny $\lambda'\:
\E' \to \E$.  Then, $\T'$, like~$\T$, is cyclic and Galois stable,
with $|\T'| = |\T|$.  Let $\E(k)$ denote the group of 
$k$-rational points of $\E$, and let 
$\beta \: \E(k) \to 
H^1(k, \T') $ be the connecting homomorphism induced on cohomology 
by the 
isogeny $\lambda '$. Then, $\T \otimes \T'\cong \mu_n$ (the group 
of $n$-th roots of unity in $k_s$).  So, we have a cup product 
pairing 
${\cupp\: H^1(k, \T) \times H^1(k, \T') \to H^2(k, \mu_n) \cong
\ _n\Br(k)}$ (the $n$-torsion of the Brauer group $\Br(k)$\,).  
Recall that the Weil-Ch\^atelet group of $\E$ 
is $H^1(k, \E)$ (which denotes, by definition, $H^1(G_k,\E(k_s)$).
This group classifies the principal homogeneous spaces of~$\E$
(see [Si, pp.~290-291]).   
Take any $\gamma\in H^1(k,\T)$, and let
$\C$ be the smooth projective curve of genus~$1$ over $k$ which 
is the principal 
homogeneous space of $\E$ determined by the image of $\gamma$ under
the canonical map $H^1(k, \T) \to H^1(k, \E)$. Then 
[CK, Th.~2.6.5] says that $\Br(k(\C)/k)$ 
coincides with the image of the map $\E(k) \to \Br(k)$ given by
$$
P \,\mapsto \, \beta(P) \,\cupp \,\gamma.
\tag 1.2
$$

This theorem will be applied here in the case where the elliptic 
curve $\E$ is 
the projective closure of $Y^2 = X^3 + c$ for $c\in k^*$.  The 
point at infinity 
$\O$  is taken as the identity element for the group operation 
on $\E$.
Then $\E$ has a subgroup $\T = \{(0,\sqrt c), (0, -\sqrt c), \O\}$, 
which is 
clearly closed under the action of $G_k$.  It takes some work to 
identify the curves 
to which the theorem applies, and to make explicit the maps used 
in the theorem.
In \S2 we prove that when 
$-3c\in k^{*2}$, the cyclic twists of $\E$ determined by $\T$
are exactly the curves $\C_f$ where $f$ is 
a diagonal form, i.e., $B = C = 0$ in (1.1) above. In~\S 3 we show 
that the desired  
isogenous curve $\E'$ is the projective closure of $Y^{\prime 2} = 
X^{\prime 3} + d$, where $d = -27c$, and we determine an isogeny 
$\lambda\:\E \to \E'$  with kernel $\T$, and also determine the dual 
isogeny $\lambda'\:\E'\to \E$, whose kernel is shown to be 
$\T' = \{(0,\sqrt d), (0, -\sqrt d), \O'\}$.  In addition, we 
determine the 
connecting homomorphism $\E(k) \to H^1(k,\T')$ induced by the 
isogeny~$\lambda'$.  In~\S 4 we apply the machinery to compute
$\Br(k(\C_f)/k)$, again in the case where $-3c\in k^{*2}$.
In \S 5 we present explicit calculations for examples with $k = \qq$, 
the rational numbers, or $k = \qq(\omega)$, where $\omega$~is a 
primitive cube root of unity.
In \S 6 we consider cyclic twists when $-3c$ is not a square
(the \lq\lq nondiagonal case"), and  show that they are
precisely the curves $\C_f$ as above where the discriminant of $f$
is not a square.  We again determine the  relative Brauer group as 
cup products as in  (1.2), though this no longer gives a 
representation of the algebras as cyclic algebras.  We give some 
explicit calculations in the nondiagonal case over $\qq$.  Finally, 
in \S 7 we turn to generalized Clifford algebras and their 
rings of quotients; we show aspects of their structure that
manifest the specialization results proved earlier.

As indicated in the paragraphs above,  we will be working with first 
cohomology groups of  Galois modules of order $3$. It is useful to 
view these modules as twisted versions of more familiar  modules, 
as follows: Let~$\A$ be a cyclic group 
of order $3$ which is a discrete $G_k$-module.  For any $b\in k^* \setminus
k^{*2}$, we define the twisted $G_k$-module $\A(b)$ to be the group 
$\A$, but with new $G_k$-action, denoted $*$, given by, 
for $a\in \A$ and $\sigma \in G_k$,
$$
\sigma*a  \ = \  \cases 
  \,\,\sigma\cdot a, &\text{if } \  \sigma(\sqrt b) = \sqrt b;\\
   (\sigma\cdot a)^{-1}, &\text{if }  \ \sigma(\sqrt b) = -\sqrt b.
\endcases
\tag 1.3
$$
If $b \in k^{*2}$ we set $\A(b) = \A$ with unchanged $G_k$-action.  
Then, for all $b,b' \in k^*$, we have ${\A(b)(b') \cong \A(bb')}$
as $G_k$-modules.  For example, take $\A = \zz_3$ (which for us 
always means 
$\zz/3\zz$ with trivial $G_k$-action).  
Then, $\zz_3(-3) \cong \mu_3$, 
the group of all cube roots of unity in $k_s$.
Also, for $\T = \{(0,\sqrt c), (0, -\sqrt c), \O\}$ as above, we have 
$\T \cong \zz_3(c) \cong \mu_3(-3c)$ as $G_k$-modules.
 Such twisted modules and 
their cohomology groups are discussed in more generality in 
[HKRT, \S 5] and [KMRT, pp.~416--418].  The following is the key 
property we need.  It is part of [HKRT, Prop.~21].

\proclaim{Proposition 1.1} Suppose $b\in k^*\setminus k^{*2}$.  
Let $L = k(\sqrt b)$, 
and let $\tau$ be the nonidentity\break $k$@-automorphism
of $L$.  The restriction map
$\res\: H^1(k,\A(b)) \to H^1(L, \A)$ is injective, and 
its image is 
$\{\delta\in H^1(L, \A)\,|\ \tau(\delta) = -\delta\}$.
\endproclaim

\bigskip\bigskip


\noindent {\bf \S2. \ Twists of $Y^2=X^3+c$ determined by $H^1(k,\T)$, diagonal case}            \medskip


	Let $k$ be a field with $\char(k)\ne 2, 3$ and let $k_s$ denote the 
separable closure of $k$.  Fix $c\in k^*$ and let $\Delta=-c/27$.  
Let $\E$ be the elliptic curve  
$$
\E\: Y^2 \, = \, X^3+c.  \tag 2.1
$$  
(More accurately, $\E$ denotes the nonsingular projective elliptic 
curve $Y^2Z=X^3+cZ^3$, which consists of the points on the affine 
curve $Y^2=X^3+c$ together with the point $\Cal O$ at infinity.
Throughout the paper, the term \lq\lq elliptic curve" always means 
a projective elliptic curve, and we routinely specify the curve
by giving an affine model, as in (2.1).) We 
denote by $\E(k)$ the group of $k$-rational points and by $\oplus$ 
the group operation on $\E$, which is defined so that $\Cal O$ is 
the identity element.  For any point $P=(r,s)\in \E(k_s)$, its 
additive  inverse $-P$ is $(r,-s)$.                                            

	Fix once and for all a choice of square root of $c$ in $k_s$, 
denoted $\sc$.  Then, $(0,\sc)$ is a point of order 3 on $\E$, 
as one can check by noting that it is an inflection point of $\E$, or 
by using the addition formula:   for $(r,s)\in \E(k_s)$,  
$$ 
(r,s)\oplus (0,\sc) \ = \, 
\left\{\aligned \textstyle \big(\frac{\,-2\sc \,r\,}{s+\sc}, 
\frac{\sc\,(s-3\sc)}{s+\sc}\big) \ \ &\text{if} \ \, s\ne -\sc,  \\ 
  \O  \hskip .65 in &\text{if} \ \, s= -\sc.  \endaligned  \right.   
\tag 2.2 
$$   
This formula can be verified for $s\ne -\sqrt c$ by 
checking that if we put $(r_1,s_1)= \big(\frac{\,-2\sc \,r\,}{s+\sc}, 
\frac{\sc\,(s-3\sc)}{s+\sc}\big)$, 
then $(r_1,s_1)\in \E(k_s)$ and the points $(r,s)$, $(0,\sc)$, 
and $(r_1,-s_1)$ are collinear.           

	Since $(0,\sc)$ has order 3 in $\E$, we have the cyclic subgroup 
of $\E(k_s)$, 
$$
\T\,=\,\{\,\Cal O, (0,\sc), (0,-\sc)\,\}.  
\tag 2.3
$$  
This $\T$ is clearly setwise invariant under the action of the 
absolute Galois group $G_k$ of $k$, and $\T$~is stabilized by 
$G_{k(\sc)}$.  As a $G_k$-module, $\T$ can be viewed as a twisted 
version of $\zz_3$ or of the group~$\mu_3$ of third roots of unity 
in $k_s$ (see (1.3) above):
	$$
\textstyle \T\, \cong \, \zz_3(c)\, \cong\,  \mu_3(\D)
 \ \ \ \text{where } \ \D= -c/27.   
\tag 2.4
$$  
Hence, we have maps 
$$
H^1(k,\T)\rA{\simeq}H^1(k,\mu_3(\D))\rA{\fres}H^1(k(\sd),\mu_3)
\rA{\simeq}k(\sd)^*\big/k(\sd)^{*3}. 
\tag 2.5
$$    
	Thus, when $\sd\in k$, $H^1(k,\T)\cong k^*/k^{*3}$.  But, when 
$\sd\not\in k$, the restriction map $\res$ is injective since 
$[k(\sd)\:\!k]$ is prime to the exponent of $\T$, and 
Proposition 1.1 yields 
	$$
H^1(k,\T)\, \cong \, \{\,z\in k(\sd)^*\big/k(\sd)^{*3}\,\big|
\,\tau(z)=z^{-1}\} ,   
\tag 2.6
$$  
where $\tau$ is the nonidentity $k$-automorphism of $k(\sd)$.

	The Ciperiani-Krashen results (cf.~[CK, \S2.6]) apply to 
projective curves 
$\C$ which are twists of~$\E$ by cohomology classes in the image 
of the map 
$$
\Psi\:H^1(k,\T)\ra H^1(k,\E) 
\tag 2.7
$$  
induced by the inclusion $\T \hookrightarrow \E(k_s)$.  We now 
determine 
those curves $\C$ in the diagonal case, where $\sd\in k$.  
(See \S 6 below 
for the case $\sd\not\in k.$) \ Fix a primitive cube root of 
unity $\w\in k_s$, 
and fix the choice of $\sd$ so that $\sc=3(2\w+1)\sd$.  
Let $P=(0,\sc)\in \T$.  
Since we are assuming $\sd\in k^*$, we have a $G_k$-module 
isomorphism  
$\T\rA{\simeq}\mu_3$ given by $P\mapsto \w$.  Take any 
$t\in k^*\setminus k^{*3}$, and choose some $\e\in k_s$ with 
$\e^3=t$.  Let 
$\g_t\:G_k\ra \T$ be given by
$$
\g_t(\r) \ = \left\{\aligned   \ \O  \ \ \  &\text{if} \ \, 
\r(\e)=\e,  \\
        P  \ \ \ &\text{if} \ \, \r(\e)=\w\e,  \\
       -P  \ \ &\text{if} \ \, \r(\e)=\w^2\e. \endaligned  \right.
\tag 2.8
$$             
	Thus, $\g_t\in Z^1(k,\T)$, and its cohomology class $[\g_t]$ in 
$H^1(k,\T)$ corresponds to $\ol t =tk^{*3}$ under the isomorphism 
$H^1(k,\T)\cong k^*/k^{*3}$.      
 
\proclaim{Proposition 2.1}  
	Assume that $\sd\in k$.  For any $t\in k^*\setminus k^{*3}$, the 
curve 
$\C_t$ corresponding to $\Psi[\gamma_t]\in H^1(k,\E)$ is  
$$
\C_t\:X^3-tY^3\,= \, -54\sd\,t^2Z^3.$$
\endproclaim                             

\demo\nofrills{\smc {Proof.}\usualspace} \         
Assume temporarily that $\mu_3\subseteq k$ (so $\sc\in k$).  
Let $K=k(\e)$ where $\e^3=t$.  Then $K$ is Galois over $k$.  
Let $G=\Gal(K/k)$ and let $\s$ be the generator of $G$ with 
$\s(\e)=\omega\e$.  
The cocycle $\g_t$ defined above is the inflation to the 
absolute Galois group $G_k$ 
of a cocycle in $Z^1(G,\T)$, which we also denote $\g_t$.  
Therefore $\g_t(\s)=P=(0,\sc)$.  
The function field $K(\E)$ of $\E$ over $K$ is $K(\rf,\sf)$ with 
$\sf^2=\rf^3+c$.  Let $*$ 
denote the twisted action of $G$ on $K(\E)$ determined by $\g_t$.  
Thus (cf.\;[Si, p.\,286]), for $f\in K(\E)$ and $\r\in G$, we have 
$\r*f=\r(f)\oo\tau_{\g_t(\r)}$, where $\r(f)$ denotes the usual 
action of $\r$ on $f$ and $\tau_{\g_t(\r)}\:\E \to \E$ is 
translation by $\g_t(\r)$, mapping 
$(r,s)\mapsto (r,s)\oplus \g_t(\r)$.  Since $\g_t(\s)=P=(0,\sc)$, 
$\tau_{\g_t(\s)}$ is given by formula (2.2) above.  
	Hence, the twisted action of $\s$ on $K(\E)$ is given by  
$$
\s* \rf \, =\, \frac{-2\sc \,\rf}{\sf+\sc}, \ \  \ 
\s*\sf \, = \, \frac{\sc\,(\sf-3\sc)}{\sf+\sc}, \ \ \
\text{and} \ \  \ \s*\e\, = \, \w\e.
$$  
This yields by straightforward calculation 
(as $9\sqrt\D=-(2\w+1)\sc$\ )   
$$
\s*(\sf-9\sd) \, = \,  \frac{-2\w^2\sc}{\sf+\sc}\,(\sf-9\sd) 
\ \  \ \text{and} \ \  \  
\s*(\sf+9\sd) \, = \,  \frac{-2\w\sc}{\sf+\sc}\,(\sf+9\sd).
$$                                                                            
	Hence, if we let $M=\e^2(\sf-9\sd)/\rf$ and $N=\e(\sf+9\sd)/\rf$, 
we find that $\s*M=M$ and $\s*N=N$, and (as $27\D=-c$)  
$$
M^3-tN^3 \, = \,  \frac{t^2}{\rf^3}\big[(\sf-9\sd)^3-(\sf+9\sd)^3\big]  
\, =  \, 
\frac{-t^2}{\rf^3}\big[54\sd(\sf^2+27\D)\big] \, = \,  -54\sd\,\,t^2.  
\tag 2.9
$$	                                                                             
	The formulas defining $M$ and $N$ show that $\rf,\sf\in K(M,N)$.  
Hence,  
$$
[K(\rf,\sf)\:\!k(M,N)]\, = \, [K(M,N)\:\!k(M,N)]\, \le \,  
[K\:\!k] \, = \, 3 \, = \, |G|.
$$      
	Therefore, $k(M,N)$ is the full fixed field of $K(\rf,\sf)$ for the 
twisted $G$-action; hence, $k(\C_t)=k(M,N)$.  Since the smooth 
projective curve $\C_t$ is determined up to isomorphism by its 
function field, we have $\C_t$ is the projective variety 
$X^3-tY^3=-54\sd\,t^2Z^3$.  This completes the proof for the 
case where $\mu_3\subseteq k$.                                  
 
	Now, suppose $\mu_3\not\subseteq k$ and let $\DD_t$ be the 
projective curve 
$X^3-tY^3=-54\sd\,t^2Z^3$ over $k$.  By [An, (3.4), (3.5), (3.8)], 
$\DD_t$ 
has Jacobian $\E$, so it is the curve determined by some 
$\th\in H^1(k,\E)$.  
$\DD_t$ clearly has an $F$-rational point in some field $F$ with 
$[F\:\!k]\,\big|\,3$, so $\res_{k\to F}(\th)=0$ in $H^1(F,\E)$.  
The 
composition $\cor\oo\res\: H^1(k,\E) \to H^1(F,\E) \to H^1(k,\E)$ is 
multiplication by 3, showing that $3\th=0$.  This tells us that 
$\th\in {}_3H^1(k,\E)$, the 3-torsion subgroup of $H^1(k,\E)$.  
Since $[k(\mu_3)\:\!k]=2$, the corestriction-restriction argument 
shows that the horizontal restriction maps in the following 
commutative diagram of 3-torsion groups are injective:  
$$
\CD
      H^1(k,\T) @>{\fres}>> H^1(k(\mu_3),\T)  \\
 	 @V{\Psi}VV         \ \  @V{\Psi_{k(\mu_3)}}VV             \\
      {}_3H^1(k,\E) @>{\fres}>>  {}_3H^1(k(\mu_3),\T) 
\endCD  
$$
	Since the earlier argument showed that 
${\Psi_{k(\mu_3)}\oo\,\res[\g_t]=\res(\th)}$, the injectivity 
implies that ${\Psi[\g_t]=\th}$, i.e., $\C_t\cong \DD_t$ over $k$.
\qed\enddemo   

\proclaim{Corollary 2.2}  
	For any $a,b \in k^*$, let $\C$ be the projective curve 
$\ \C\:Z^3 = aX^3+bY^3$.  
Then $\C$ is a homogeneous space for the elliptic curve with affine 
model $\E\:Y^2=X^3+c$ where $c=-\frac{27}4a^2b^2$,  and $\C$ lies 
in the image of $\Psi\:H^1(k,\T)\to H^1(k,\E)$.  Specifically, 
$\C$ is $\Psi[\g_t]$ for $t=-b/a$ and also for $t=a$ and for $t=b$.
\endproclaim

\demo\nofrills{\smc {Proof.}\usualspace} \   We have $\D=-c/27 = 
(ab/2)^2\in k^{*2}$.  Take $\sd=ab/2$.  Then  Proposition 2.1 with 
${t=-b/a}$ shows that $\C_{-b/a}$ is the curve 
$X^3-(-b/a)Y^3=-54(ab/2)
(-b/a)^2Z^3$, that is, ${aX^3+bY^3=(-3bZ)^3}$, which is clearly 
isomorphic to $\C$.  (If we made the other choice of $\sd$, then 
$\C_{-b/a}$ is ${aX^3+bY^3=(3bZ)^3}$, which is also isomorphic 
to $\C$.) \;  Likewise, $\C_a$ can be rewritten as 
$X^3=aY^3+b(-3aZ)^3$ and 
$\C_b$ as $X^3=a(-3bZ)^3+bY^3$, so $\C\cong \C_a$ and $\C\cong \C_b$.   
\qed\enddemo

\bigskip\bigskip                                   

\vfill\eject       


\noindent {\bf \S3. \ Dual isogenies and a connecting homomorphism}            \medskip


	For the elliptic curve 
$\E\:y^2=x^3+c$  and its subgroup 
$\T=\big\{\,\O, \; (0,\pm\sc)\,\big\}$ 
as in \S 2, there is, up to isomorphism, a unique elliptic curve 
$\E'$ and isogeny $\l\:\E\to \E'$ with $\ker(\l)=\T$.  That is, 
we have a short exact sequence 
$$
0\ra \T \ra \E \rA{\l}  \E' \ra 0.   
\tag 3.1
$$
There is further an isogeny dual to (3.1) \, (cf.~[Si, Ch.~III, \S6]) 
$$
0\ra \T' \ra \E' \rA{\l'} \E \ra 0,   
\tag 3.2
$$  
and we need to determine the associated connecting homomorphism 
$H^0(k,\E) \to H^1(k,\T')$.  First, we obtain explicit descriptions 
of $\l$, $\E'$, $\T'$, and $\l'$.

	The paper [V] by V\'elu describes how to compute the isogeny
 with kernel $U$ for any finite subgroup $U$ of an elliptic curve.  
When $U$ is our group $\T$, his prescription is as follows:  For 
a point~$P$ of $\E(k_s)$, write $P=(r(P),s(P))$; then for 
$P\not\in \T$,
$$
\l(P) \ = \ \bigg(r(P)+\!\!\!\sgm{Q\in \T\setminus\{\O\}}\!\!\!
\big(r(P\oplus Q)-r(Q)\big), 
\ \ \ s(P)+\!\!\!\sgm{Q\in \T\setminus\{\O\}}
\!\!\!\big(s(P\oplus Q)-s(Q)\big)\bigg).    
\tag 3.3
$$  
	If $P\in \T$, then $\l(P)=\O'$, the point at infinity, which is the 
identity element of $\E'$.  Since ${\T\setminus\{\O\}=
\big\{(0, \sc),(0,-\sc)\big\}}$, for $(r,s)\in \E(k_s)$, we set 
$(r_1, s_1)=(r,s)\oplus (0,\sc)$ and\break $(r_2,s_2)=(r,s)\oplus (0,-\sc)$.  
	Then, by (2.2) above and its analogue for $(0,-\sc)$, 
$$
\l(r,s) \, =\,  (r+r_1+r_2,\; s+s_1+s_2) \, = \,  
\left(\frac{r^3+4c}{r^2}, \ \frac{s^3-9cs}{r^3}\right)  \ \ \ 
\text{if} \ (r,s)\not\in \T,  
\tag 3.4
$$ 
and $\l(r,s)=\O'$ if $(r,s)\in \T$.  
	(This formula for the isogeny $\l$ also appears as an exercise in 
[Ca,~p.~65, \#5].)
	By [V, (11)] the elliptic curve $\E'$ has Weierstrass model  
$$
\E'\:\,Y^2\,=X^3+d \ \ \ \text{where} \ d=-27c.
$$
	(One can  verify this formula for $\E'$ by checking that all 
the points $\l(r,s)$ lie on this curve.)

	The isogeny $\l$ has degree $|\T|=3$.  Therefore the kernel $\T'$ 
of the dual isogeny $\l'$ is the image under $\l$ of the points  
of $\E$ of order 3 (cf.~[Si, p.~84, Theorem 6.1(a)]).  The order 
$3$ points $(\rh,\sh)$ are obtainable 
as the inflection points of $\E$, i.e., points where 
$d^2y\big/dx^2=0$.  
An easy calculus exercise shows that these are the points with 
($\rh=0$ and $\sh^2=c$) or ($\rh^3=-4c$ and $\sh^2=-3c$).  Hence, 
by computing $\l(\rh,\sh)$ we find  
$$
\T' \,= \,  \big\{\, \O', \ (0,\sqrt{d}), \ (0,-\sqrt{d}) \,\big\}.   
\tag 3.5
$$
	(We could alternatively have found the points $(\rh,\sh)$ by using the 
triplication formula for $\E$ given in [Si, p.~105, Ex.~3.7(d)].)  By 
repeating the process used to find $\l$, with $d$ replacing $c$, we 
obtain an isogeny $\wt{\l}\:\E'\to \E''$  with $\ker(\wt{\l})=\T'$ 
where $\E''$ is the elliptic curve $\E''\:Y^2=X^3-27d=X^3+27^2c$.
	To get back to $\E$, use the isomorphism $\iota\:\E''\to \E$ 
given by $(r,s)\mapsto (r/9,\, s/27)$ and $\O''\mapsto \O$.  Thus, 
the isogeny $\l'=\iota\oo\wt{\l}\: \E'\to \E$ with kernel $\T'$ is 
given by 
$$
\l'(r',s') \,= \, \left(\frac{r'{}^3+4d}{9r'{}^2}, \ 
\frac{s'{}^3-9ds'}{27r'{}^3}\right) \ \ \ \text{if} \ (r',s')\not\in \T',  
\tag 3.6
$$  
and $\l'(r',s')=\O$ if $(r',s')\in \T'$.                                
  
\def\t{\widehat t}

\proclaim{Lemma 3.1}  
	Let $(r,s)\in \E(k_s)$ with $s\ne -\sc$.  Choose $t\in k_s^*$ 
with $t^3=s+\sc$,  and set $\t=r/t$.  Let $r'=6\sc\big/(t-\t\,)$ and 
$s'=9\sc(t+\t\,)\big/(t-\t\,)$.  Then $(r',s')\in \E'(k_s)$ and 
$\l'(r',s')=(r,s)$.
\endproclaim 

\demo\nofrills{\smc {Proof.}\usualspace} \         
	Since $s\ne -\sc$, it follows that $t\ne 0$, and so $\t$ is 
well-defined.  Also, by hypothesis we 
have  $\t\,^3=r^3/t^3=s-\sc \ne s+\sc=t^3$.  Therefore, $\t\ne t$, 
assuring that 
$r'$ and $s'$ are well-defined.
	(The formulas for $r'$ and $s'$ were found by noting that if 
$(\wt{r}',\wt{s}')\in \E'(k_s)$ and 
${\l'(\wt{r}',\wt{s}')=(\wt{r},\wt{s})}$, 
then, from (3.6),  
$\wt{s}+\sc=$ $[(\wt{s}'{}^3-9d\wt{s}')\big/ {27\wt{r}'{}^3}]+\sc=
\big[(\wt{s}'+9\sc)\big/3\wt{r}'\big]^3$, and likewise
${\wt{s}-\sc=\big[(\wt{s}'-9\sc)\big/3\wt{r}'\big]^3}$.  These 
equations can 
be solved to recover $\wt{r}'$ and $\wt{s}'$.)  Note that as 
$(t-\t\,)^3=2\sc-3(t-\t\,)r=-\big[3(t+\t\,)^2(t-\t\,)-8\sc\big]$, 
we obtain $s'{}^2-r'{}^3=-27c=d$; so $(r',s')\in \E'(k_s)$.   
Likewise, straightforward calculations show that 
$(r'{}^3+4d)\big/9r'{}^2=r$ and 
$s'(s'{}^2-9d)\big/27r'{}^3=s$; so, $\l'(r',s')=(r,s)$.   \qed\enddemo

\medskip

	Since $\sqrt{d}=3\sqrt{-3}\sc$, it follows that $\T'$ is a twisted 
$G_k$-module, and 
$$
\T'\, \cong \,  \mu_3(c).  \tag 3.7$$   Hence, we have 
maps  
$$
H^1(k,\T') \rA{\simeq} H^1(k,\mu_3(c)) \rA{\fres} H^1(k(\sc),\mu_3) 
\rA{\simeq} k(\sc)^*/k(\sc)^{*3}.  
\tag 3.8
$$      
	Define the homomorphism 
$$
\a\:\E(k) \to k(\sc)^*/k(\sc)^{*3},  
\tag 3.9
$$ 
to be the composition of maps  
$$
\E(k) \rA{\partial} H^1(k,\T') \ra k(\sc)^*/k(\sc)^{*3}, 
\tag 3.10
$$ 
 where $\partial$ is the connecting homomorphism $\E(k)=H^0(k,\E)
\to H^1(k,\T')$ arising from  the short exact sequence (3.2), and the 
second map is the monomorphism which is the composition of maps in (3.8).

\proclaim{Proposition 3.2}  
	For $(r,s)\in \E(k)$ with $s\ne -\sc$, one has 
$$
\a(r,s) \, = \, \,\ol{s+\sc\,}  \, \in  \, k(\sc)^*/k(\sc)^{*3}.
$$  
Also, if $\sc\in k$, then $\a(0,-\sc)=\ol{4c}$.
\endproclaim

\demo\nofrills{\smc {Proof.}\usualspace} \         
To be consistent with choices of square roots of $c$ and $d$ and  
cube roots of unity, let~$\w$ be a fixed cube root of unity, 
and let $\eta = 2\w +1$, which is a fixed square root of 
$-3$.  Then set $\sqrt{d}=3\eta\sc$.  
We use the $G_k$-module isomorphism $\T'\to \mu_3(c)$ given by 
$(0,\sqrt d)\mapsto \w$.  Since connecting homomorphisms are 
compatible with restriction maps and the map  
${\res\: H^1(k,\T') \to H^1(k(\sc), \T')}$ is injective, it suffices 
to prove the proposition when $\sc\in k$.  Assume this.  For any 
$u \in k^*$ and any fixed $\cbrt{u}\in k_s^*$, we have a cocycle 
$\g_u\in Z^1(k,\T')$ given by, for $\s\in G_k$, 
	$$
\g_u(\s) \, =\  \left\{\,\aligned  \O' \ \ \ \ \ \ \ &\text{if} \ \, 
\s(\cbrt{u})= \cbrt{u},   \\ 
(0,\sqrt{d}) \ \ \ \,  &\text{if} \ \, \s(\cbrt{u})=\w\,\,\cbrt{u},   \\  
(0,-\sqrt{d}) \ \  &\text{if} \ \, \s(\cbrt{u})=\w^2\,\,\cbrt{u}.
\endaligned 
\right. 
$$     
	Our isomorphism $H^1(k,\T')\to k^*/k^{*3}$ of (3.7) is given by 
$[\g_u] \mapsto \ol u$.  Thus, it suffices to prove that the 
connecting homomorphism $\partial\: \E(k) \to H^1(k,\T')$ maps 
$(r,s)\mapsto [\g_{s+\sc}]$ when $s\ne -\sc$, and 
$\partial(0,-\sc)=[\g_{4c}]$.

	For $(r,s)\in \E(k)$ with $s\ne -\sc$, choose $t\in k_s^*$ with 
$t^3=s+\sc$, and let $\t= r/t$, ${\rho=t+\t}$, and $\d=t-\t\ne 0$.  
Set $(r',s')=(6\sc/\d, \, 9\sc \rho/\d)$.  By Lemma 3.1, 
$(r',s')\in \E'(k_s)$ and ${\l'(r',s')=(r,s)}$.  Then 
$\partial (r,s)$ is the cohomology class of the cocycle 
$\zeta\:G_k\to \T'$ given by $\s\mapsto \s(r',s')\ominus (r',s')$, 
where $\ominus$ denotes subtraction in $\E'$.  Take any $\s\in G_k$.  
	If $\s(t)=t$, then also $\s(\t)=\t$, as $\t= r/t$ and $x\in k$, so 
$\s(r',s')=(r',s')$.  Hence, $\zeta(\s)=\O'=\g_{s+\sc}(\s)$. 
	Suppose next that ${\s(t)=\w t=\frac12(-1+\eta)t}$.  Then, as 
$t\,\t=r$, we have $\s(\t)=\w^{-1}\t=\frac12(-1-\eta)\t$, so that 
${\s(\rho)=\frac12(-\rho+\eta\d)}$ and $\s(\d)=\frac12(-\d+\eta \rho)$, 
which is nonzero as $\d\ne 0$.  Hence,  
	$$
\s(r',s')\,= \, \s\big(6\sc/\d, \ 9\sc \rho/\d\big) 
\, = \,  \big(12\sc\big/(-\d+\eta \rho\,), \ 
9\sc(-\rho+\eta\d)\big/(-\d+\eta \rho)\big).
\tag 3.11
$$   
Note that $s'+\sqrt{d}=3\eta\sc(-\eta \rho+\d)\big/\d\ne 0$.           

	The analogue to formula (2.2) (replacing $c$ by $d$) combined 
with (3.11) yields
$$
\aligned 
(r',s')\oplus (0,\sqrt{d}) \ 
	   & = \ \big(\textstyle\frac{-2\sqrt d r'}{\,s'+\sqrt d\,},  
\frac{\sqrt d(s'-3\sqrt d)}{s'+\sqrt d}\big)  \\
	    & = \ \big(\textstyle\frac{-36\eta c}{\,9\sc \rho\,+\,3\eta\sc\d\,},  
\frac{\,3\eta\sc(9\sc \rho\,-\,9\eta\sc\d)\,}{9\sc \rho\,+\,3\eta\sc\d}\big) \\
		&	= \ \big(\textstyle\frac{12\sc}{\,\eta \rho-\d\,}, 
\frac{\,9\sc(\rho-\eta\d)\,}{-\eta \rho+\d}\big)  \\
		&	= \ \s(r',s'). \endaligned 
$$
	Thus, we have $\s(r',s')\ominus (r',s')=(0,\sqrt d)$ when 
$\s(t)=\w t$.
	Likewise, if $\s(t)=\w^2t$, then the analogous calculation 
(with $-\sqrt d$ and $-s$ replacing $\sqrt d$ and $s$) shows that 
${\s(r',s')\ominus(r',s')=(0,-\sqrt d)}$.  So our cocycle 
$\zeta=\g_{\,t^3}=\g_{y+\sc}$, showing that 
$\partial(r,s)=[\g_{s+\sc}]$ 
whenever $s\ne-\sc$.  Since\break $\partial$ is a homomorphism with 
$\partial(0,\sc)=[\g_{2\sc}]$ and $(0,-\sc)=(0,\sc)\oplus(0,\sc)$,
 we have\break ${\partial(0,-\sc)=[\g_{(2\sc)^2}]=[\g_{4c}]}$.  
\qed\enddemo

\bigskip\bigskip


\noindent {\bf \S4. \ Relative Brauer Groups $\Br(k(\C)/k)$ 
for $\C\:Z^3=aX^3+bY^3$}            \medskip


For any $a, b\in k^*$, let $\C$ be the smooth projective genus $1$ 
curve over $k$, 
$$
\C\: \, Z^3 \, =\,  a\,X^3 \, + \, b\,Y^3\,,\,
\tag 4.1
$$ 
and let $k(\C)$ be the function field of $\C$ over $k$.  We can now 
describe the relative Brauer group $\Br(K(\C)/k)$, as the target 
of a homomorphism from the $k$-rational points of the Jacobian of 
$\C$, which is the elliptic curve $\E$ with Weierstrass model 
$$
\E\: \, Y^3 \,= \,  X^3 \,+\,  c\,,\quad \text{where} 
\ \  c \,=\, \textstyle- \frac{27}4 a^2 b^2\,.
\tag 4.2
$$
The algebras in $\Br(k(\C)/k)$ are symbol algebras  when 
$\mu_3\subseteq k$, 
and are cyclic algebras otherwise.  
Throughout this section, let $\w$ denote a 
fixed primitive cube root of unity in $k_s^*$.
If $\w \in k^*$, let 
$(t,u;k)_{\w}$, or 
simply $(t,u)_{\w}$, denote the $9$-dimensional symbol algebra 
over $k$ with generators $i$, $j$ and relations $i^3=t$, $j^3=u$, 
and $ij=\w ji$.

\proclaim{Proposition 4.1}  
Let $k$ be a field with $\mu_3 \subseteq k$, and let    $a,b\in k^*$. 
For the projective curve  $\C\:Z^3=aX^3+bY^3$  over $k$ 
with Jacobian  $\E$ as in $(4.2)$, 
 there is a surjective group homomorphism  
	$$
\varphi\:\E(k)\to\Br(k(\C)/k) \ \ \text{given by}  
\ \, (r,s) \mapsto [(a,s+\sqrt{c\,})_{\w}]\,. \tag 4.3
$$   
In other words, 
$$
\Br(k(\C)/k)  \, =  \,\{\,[(a,s+\sqrt{c\,})_{\w}] \; 
\big| \ \ol{s+\sqrt{c\,}\,}\in \im(\a) \subseteq k^*\big/k^{*3}\,\},
$$ 
where $\a$ is the map in $(3.9)$ and Proposition~$3.2$. 
\endproclaim

\demo{Proof}  (When $k$ is perfect, this was proved  in 
[H$_1$, Cor.~ 1.2]
and [H$_2$, Th.~ 1.2 and Cor.~ 2.2] 
using the generalized Clifford algebra of the binary cubic form $aX^3 + bY^3$.  
The case when $k$ is not perfect is deducible from the 
perfect case by working back from the perfect hull of $k$.  We give
an alternative approach here based on the theorem of 
Ciperiani and Krashen, since this approach is more readily applicable
in other cases.  See Proposition 4.4 and \S 6 below.)  

Note that $\sqrt c\in k$, as $\w \in k$.  Let $\T$ be the cyclic 
subgroup  $\{\O, \,(0, \pm \sqrt c)\}$ 
of $\E(k)$.  By Corollary~2.2 above, $\C$ is the 
homogeneous space of $\E$ corresponding to $\Psi[\gamma_a]$ in 
$H^1(k, \E)$, where $\gamma_a \in Z^1(k, \T)$ is as in (2.8).  
We have the isogeny (3.1)  with kernel $\T$, and its dual 
isogeny (3.2), with kernel $\T' = \{ \O', (\pm \sqrt d)\}$, where 
$d = -27c$. Let $\partial\: \E(k) \to H^1(k, \T')$ be 
the connecting homomorphism arising from the dual isogeny.
By [CK, Th.~ 2.6.5], $\Br(k(\C)/k)$ coincides with the image of the map 
$$
\beta\: \E(k) \to \Br(k) \ \ \ \ \text{given by} \ \ \ \ 
\beta(r,s) \, =\,  \partial (r,s) \cupp [\gamma_a]\,,
\tag 4.4
$$ 
where the cup product $\cupp$ maps $H^1(k, \T') \times H^1(k, \T)\to 
H^2(k, \T'\otimes \T) \cong H^2(k, \mu_3) \cong {}_3\Br(k)$.
Now, $\T'\cong \zz_3(c) \cong \mu_3$ as $G_k$-modules, since 
$\sqrt c\in k$ and $\mu_3 \subseteq k$.  Hence, $H^1(k,\T')\cong 
k^*/k^{*3}$ 
and, as shown in Proposition 3.2, under this isomorphism the map 
$\partial$ corresponds to $\alpha\:\E(k) \to k^*/k^{*3}$ given 
by $(r,s) \mapsto \overline{s+\sqrt c \,}$.
Also, $\T\cong \zz_3(d) \cong \mu_3$, and in the isomorphism 
$H^1(k,\T) \cong k^*/k^{*3}$, $[\gamma_a]$~corresponds to $\ol a$.  
The cup product $H^1(k, \mu_3) \cupp H^1(k, \mu_3) 
\to {}_3\Br(k)$ gives rise  to the map 
$$
\kappa\: k^*/k^{*3} \times  k^*/k^{*3} \to {}_3\Br(k)
$$
given by 
$(\ol t, \ol u) \mapsto (u,t)_\omega$. (For, by [Se, p.~207,
 Prop.~ 5], 
the cup product maps the pair $(\ol t, \ol u)$ to the Brauer class 
of the cyclic algebra $(k(\root 3\of t)/k, \sigma, u)$, where 
$\sigma(\root 3\of t) = \omega\root 3 \of t$.  
This algebra is generated by $i$ and $j$, where $j^3 = t$, 
$i^3 = u$ and $iji^{-1} = \omega j$, so $ij = \omega ji$.)  Thus, 
for any $(r,s) \in \E(k)$, we have
$\beta(r,s) = \kappa(\ol{s+\sqrt c\,}, \ol a) = (a, s+\sqrt c)_\omega$.    
\qed \enddemo

\noindent{\bf Remark 4.2}.  (i) \; If the symbol algebra 
$(a,b)_{\w}$ is nonsplit, then Proposition 4.1 shows that 
$\Br(k(\C)/k)$ is nontrivial since $[(a,b)_{\w}]\in\Br(k(\C)/k)$.  
In fact, since $\w\in k$, we have $\sqrt c\in k$, so the point 
$(0,\sqrt c)\in\E(k)$.  This implies that 
$[(a,2\sqrt c)_{\w}]\in\Br(k(\C)/k)$.   But, we have 
$$
\textstyle(a,2\sqrt c)_{\w} \, =\, 
\big(a,2\sqrt {\botsmash{-\frac{27}4a^2b^2}}\big)_{\w}
\, \cong \, \big(a,\sqrt{-3}^3ab\big)_{\w} 
\, \sim\,  (a,a)_{\w}\otimes(a,b)_{\w} 
\,\sim\, (a,b)_{\w}.
$$
Note here that $-1$ is a cube in $k$, so 
$(a,a)_{\w}\cong(a,-1)_{\w}$ is split.  
	We  point out also that the converse of this remark is not 
necessarily true in general.   Corollary 4.6\,(i) below can be 
used to provide counterexamples.

	(ii) \; In Proposition 4.1, we can replace $a$ in the symbol algebras by $b$
or by $-b/a$ if this is more convenient for computation.  This follows
from Corollary 2.2, since the curve $\C$ is represented by 
$\Psi[\gamma _b]$ and $\Psi[\gamma_{-b/a}]$ as well as by 
$\Psi[\gamma_a]$.   

\medskip

	Now, we consider the case that $\mu_3 \not \subseteq k$.  In 
order to determine $\Br(k(\C)/k) $  for any field $k$, we~need 

\proclaim{Lemma 4.3}  
Let $k$ be a field with $\mu_3 \not \subseteq k$,  let $L=k(\w)$, and let 
$\sigma\in \Gal(L/k)$ with $\sigma \ne \id$.  
 If $T=L(\cbrt{\xi\,})$ for $\xi\in L^*\setminus L^{*3}$,  
then$:$ 
\roster
\item "{\rm (i)}" There is a subfield $S$ of $L$ such that 
$S/k$ is a cyclic extension of degree $3$ if 
and only if $\sigma(\ol{\xi})=\ol{\xi}\,^{-1}$ in $L^*/L^{*3}$.
\item "{\rm (ii)}" If $\s(\xi)=\xi^{-1}r^3$ for some $r\in k^*$, then $S=k(u)$, 
where $u$ has  minimal polynomial ${X^3-3rX-(\xi+\s(\xi))}$
over $k$.
\endroster 
\endproclaim

\demo{Proof} (i) This result is known, see Albert ([A, Th.~ 2]). 
In cohomological terms, since $\mu_3 \not\subseteq k$, we have 
$\zz_3\cong \mu_3(-3)$ 
as $G_k$-modules.  Hence, there are maps 
$$
{H^1(k,\zz_3) \rA{\simeq} H^1(k, \mu_3(-3)) 
\hookrightarrow H^1(L,\mu_3)}.
$$ 
This yields an injection $H^1(k, \zz_3)\to L^*/L^{*3}$, 
whose image is 
${\{\, \ol\xi\in L^*/L^{*3} \,|\ \s (\ol\xi) = \ol \xi\,^{-1} \}}$ 
by Proposition 1.1 . 
 The cyclic subgroups of $H^1(k, \zz_3)$ classify the cyclic field 
extensions of 
$k$ of degree $3$, and in field terms, the preceding  map takes a 
cyclic field 
extension $S$ of $k$  to $\ol\xi$ where $S\cdot L = L(\cbrt \xi)$.

  (ii) Let $\e\in T$ with $\e^3=\xi$.  The map $\s$ extends to an 
automorphism 
$\wt\s$ of $T$ of order $2$ given by $\wt\s(\e) = r\e^{-1}$, and 
the fixed field of $\wt\s$ is $S$.  Let 
$u = \e+ \wt\s(\e) = \e + r\e^{-1} \in S$.  
Then, $u\notin k$ since $\e$ satisfies a polynomial of degree 
$2$ over $k(u)$ but $[L(\e)\!:\!L] = 3$.  So, $S = k(u)$.  
By computing $u^3$, 
we see that $u$ is a root of $X^3 -3rX - (\xi + \s(\xi)) \in k[X]$.  
Since $[k(u)\!:\!k] = 3$, this polynomial must be the minimal
polynomial of $u$ over $k$.  Let $u' = \w\e + \wt\s(\w\e) = 
\w\e + \w^{-1}r\e^{-1}$. Then $u'$ is another root of $\min(u,k)$ 
in $S$.
So, there is a nontrivial $k$-automorphism of $S$ sending $u$
to $u'$.  This shows directly that $S$ is cyclic Galois over $k$.
\qed \enddemo   

	Lemma 4.3 enables us to give the analogue to  Proposition 4.1 when 
$\mu_3 \not \subseteq k$.  In that case we let 
$L = k(\omega) = k(\sc)$, and work down from $3$-Kummer 
field extensions of $L$ to cyclic extensions of $k$.

\proclaim{Proposition 4.4}  
Let $k$ be a field with $\mu_3 \not\subseteq k$.  For 
$\C\: Z^3 = aX^3 + bY^3$
and $c$, $\E$ as in $(4.2)$, there is a surjective group 
homomorphism  
$$
\varphi\:\E(k)\to\Br(k(\C)/k) \ \ \text{given by} \ \, (r,s) \mapsto 
[(S/k,\tau,a)]
$$ 
where $S$ is the  cyclic field extension of $k$ lying in 
$k\big(\w, \cbrt{s+\sqrt c\,}\big)$ described in Lemma 4.3.
\endproclaim

\demo\nofrills{\smc {Proof.}\usualspace} \ 
The result of [CK] still applies, just as for Proposition 4.1, and 
it shows that\break $\Br(k(\C)/k) = \im(\beta)$ for $\beta$ the map 
of (4.4).
But now, as $\mu_3 \not \subseteq k$, we have $\sqrt c\notin k$.  Let 
$L = k(\sqrt c) = k(\w)$, and let $\s$ be the nonidentity 
$k$-automorphism of $L$.  We have $\T'\cong \zz_3
\cong \mu_3(c)$ as $G_k$-modules,
and $\T\cong \mu_3$.  So the cup product in 
(4.4) corresponds to the cup product 
$H^1(k, \zz_3) \times H^1(k, \mu_3)\to {}_3\Br(k)$.  
Here $[\gamma_a] \in H^1(k,\T)$ corresponds
to $\ol a\in k^*/k^{*3} \cong H^1(k,\mu_3)$.   

Now, $H^1(k, \zz_3)\cong \Hom(G_k, \zz_3)$, which classifies 
cyclic degree $3$ field extensions of $k$  with a specified generator 
of the Galois group.  The restriction map sends $H^1(k, \zz_3)
\cong H^1(k, \mu_3(c))$  injectively into $H^1(L, \mu_3) \cong
L ^*/L^{*3}$.  
By Proposition 1.1, the image of $H^1(k, \zz_3)$ in $L^*/L^{*3}$ 
consists of those $\ol\xi$ with 
$\s(\xi) =r^3\xi^{-1}$ for some $r\in L^*$.  Since $\s(r^3) = r^3$, 
$r$ can be chosen so that $r\in  k^*$.    For such a $\ol \xi$ the 
field extension  of $k$ associated to the inverse image of 
$\ol\xi$ in $H^1(k,\zz_3)$ is the field~$K$ with 
$K\cdot L = L(\cbrt\xi)$, which is described in Lemma 4.3.  
For $(r,s) \in \E(k)$, by definition of the map~ $\a$ in (3.10), 
$\partial(r,s)\in H^1(k, \T')$ maps to $\alpha(r,s)
=\ol{s+\sqrt c\,}$ in $L^*/L^{*3}$ (cf.~Proposition 3.2).  So, 
the associated cyclic field extension of $k$ is the $S$ of 
Lemma 4.3 lying in 
$L(\cbrt{s+\sqrt c\,})$. Thus, by [Se, p.\,204, Prop.~ 2], 
$\partial(r,s) \cupp
[\gamma_a]$ is the Brauer class of the cyclic algebra 
$(S/k,\tau, a)$, where $\tau$ 
is the restriction to~$S$ of the $L$-automorphism of 
$L(\cbrt{s+\sqrt c\,})$ sending 
$\cbrt{s+\sqrt c\,}$ to $\omega\,\, \cbrt{s+\sqrt c\,}$. 	
\qed \enddemo 

\medskip

\noindent{\bf Remark 4.5}.  \  Lemma 4.3\,(ii) shows that 
if $(r,s)\in \E(k)$, then the corresponding field $S$ is~$k(u)$, 
where 
${u\,=\,
\cbrt{s+\sqrt{c\,}\,}\,+\,\,\cbrt{s-\sqrt{c\,}\,},}$
(with the cube roots chosen so that their product is $r$); 
the minimal polynomial of $u$ over $k$ is $\min(u,k)=  X^3-3rX-2s$. 
	Note that this minimal polynomial has discriminant 
$-4(-3r)^3-27(-2s)^2=4\dd27(r^3-s^2)=27^2a^2b^2$, which is a square 
in $k$.  Thus, the Galois group of the splitting field $S$ of 
$\min(u,k)$ over $k$ must be cyclic of order 3.

\medskip

\proclaim{Corollary 4.6}  
	Let $\C\: Z^3=aX^3+bY^3, \ a,b\in k^*$, be a curve over $k$, 
and $\E$ the Jacobian of~$\C$.    For the map 
$\a\:\E(k)\to k(\sc)^*/k(\sc)^{*3}$ in $(3.9)$, if $\im(\a)$ 
is trivial, then  $\Br(k(\C)/k)=\{\,0\,\}.$
\endproclaim

\demo\nofrills{\smc {Proof.}\usualspace} \ The definition of $\a$ 
shows that if $\a$ is trivial, then so is the connecting 
homomorphism $\partial\: \E(k) \to H^1(k, \T')$.  When this occurs, 
the epimorphism $\b\: \E(k) \to \Br(k(\C)/k)$ of (4.4) is clearly 
also trivial, and so its image must be trivial. 
\qed \enddemo

\bigskip\bigskip


\noindent {\bf \S5. \ $\Br(k(C)/k)$ for $k$ a global field with 
specific examples, 
diagonal case}            \medskip


	In this section, we study various useful facts on elliptic curves 
of the form 
$\E\:Y^2=X^3+c$, which are Jacobians of curves $\C\:Z^2=f(X,Y)$ where 
$f(X,Y)$ is a binary cubic form.  
Based on these facts, we give assorted examples of relative 
Brauer groups 
$\Br(k(\C)/k)$ where $\C$ is a diagonal cubic curve 
$Z^3 = aX^3+bY^3$ over 
$k=\qq$ or $\qq(\w)$.  
Examples of nondiagonal case will be given at the end of  \S 6.

\proclaim{Lemma 5.1}  
	Let $\l\: A \to B$ and $\l'\:B\to C$ be group homomorphisms such 
that the kernels 
and cokernels of $\l$ and $\l'$ are all finite.  Then  
$\ker\,(\l'\oo \l)$ and 
$\coker\,(\l'\oo \l)$ are finite, and 
	$$
\frac{\,\big|\coker\,(\l'\oo \l)\big|\,}{\big|\ker\,(\l'\oo \l)\big|}
\, =\,
\frac{\,\big|\coker\,(\l')\big|\dd\big|\coker\,(\l)\big|\,}
{\big|\ker\,(\l')\big|\dd\big|\ker\,(\l)\big|}.   \tag 5.1 
$$
\endproclaim                               

\demo\nofrills{\smc {Proof.}\usualspace} \ It is easy to check
that there is  is an exact 
sequence with the obvious maps (cf.\; [Mi, Lemma A.2, p.\,86]),
$$
0 \to \ker(\l) \to \ker(\l'\oo\l) \to \ker(\l') \to \coker(\l) 
\to \coker(\l'\oo\l) \to \coker(\l') \to 0.
$$  
Since the kernels and cokernels of $\l$ and $\l'$ are all 
assumed to be  finite, it 
is easy to check (or see [Mi, Lemma 3.7, p.\,79]) that 
$$
\big|\ker(\l)\big|\, \big|\ker(\l')\big| \, 
\big|\coker(\l'\oo \l)\big| \, =  \, 
\big|\ker(\l'\oo \l)\big| \,\big|\coker\,(\l)\big|
\, \big|\coker\,(\l')\big|.
$$
This yields the formula in (5.1).   \qed\enddemo 

	Let $\E\:Y^2=X^3+c$ and $\E'\: Y^2=X^3+d$, where $d=-27c\in k^*$, 
be  elliptic curves over $k$.  Recall that the 
{\it rank} of $\E(k)$, call it $\rk$,  
is the rank torsion-free of $\E(k)$, i.e., 
$\rk = \dim_{\qq}(\E(k)\otimes_{\zz}\qq)$.   
Notice  that $\E(k)$ and $\E'(k)$ have the same rank since 
they are isogenous.  
	
	In Lemma 5.1, if we put $A=C=\E(k)$ and $B=\E'(k)$, we have

\proclaim{Proposition 5.2}  
	Let $k$ be a global field.  Let $\l$, $\l'$ be the maps as in 
$(3.1)$, $(3.2)$, 
the map $\a$ as in Proposition 3.2, and the map $\a'$ analogous 
to $\a$, $d$ replacing $c$.  
If $\E(k)$ has rank $\rk$, then 
$$
3^{\rk} \, =  \, \frac{\big|\im\,(\a)\big|\dd\big|\im\,(\a')\big|}
{\,\big|\ker\,(\l')\big|\dd\big|\ker\,(\l)\big|\,}.
$$
\endproclaim

\demo\nofrills{\smc {Proof.}\usualspace} \         
	By the Mordell-Weil Theorem, one has 
$\E(k) \cong (\ds{j=1}^{\rk} \zz) 
\bigoplus \big(\ds{i=1}^\ell \zz_{{p_{{}_i}}^{\!\!n_i}})$, 
and so  
$$
\E(k)/3\E(k)  \, \cong \,  \big(\ds{j=1}^{\rk} \zz/3\zz\big) 
\bigoplus 
\big(\ds{i=1}^\ell 
\zz_{{p_{{}_i}}^{\!\!n_i}}/3\zz_{{p_{{}_i}}^{\!\!n_i}}\big).
$$
  If the prime $p_i$ is different from 3, then 
$\zz_{{p_{{}_i}}^{\!\!n_i}}/3\zz_{{p_{{}_i}}^{\!\!n_i}}=0$.  
On the other hand, since 
${\zz_{{3}^{n_i}}/3\zz_{{3}^{n_i}}\cong\zz/3\zz}$, 
it follows that 
$$
[\E(k):3\E(k)] \,  =  \, 3^{\rk+t}
$$ 
where $t$ is the number of $i$ with $p_i = 3$.
	If $\E(k)_3$ denotes the 3-torsion subgroup of $\E(k)$, then 
$|\E(k)_3| = 3^t$ and thus 
$$
[\E(k):3\E(k)]  \, = \,  3^{\rk}\!\cdot|\E(k)_3|. \tag 5.2
$$  
	Now, since $\big|\coker\,(\l'\oo \l)\big| = 
\big|\E(k)/3\E(k)\big|=
 [\E(k)\:\!3\E(k)]$ and $\big|\ker\,(\l'\oo \l)\big| = 
\big|\E(k)_3\big|$, 
it follows from (5.2) that $\frac{\big|\coker\,(\l'\circ \l)\big|}
{\big|\ker\,(\l'\circ \l)\big|}= 3^{\rk}$.  On the other hand, 
we have 
$$
\coker\,(\l') \, = \,  \E(k)/\im(\l')  \, \cong \, 
 \E(k)/\ker(\a)  \, \cong \,  \im(\a).
$$ 
and analogously $\coker\,(\l)\cong \im(\a')$.  This completes the 
proof.   \qed\enddemo

	The following corollary describes the relationship between the rank 
$\rk$ of $\E(k)$ and the images of the maps $\a$ and $\a'$ 
respectively.  
This observation leads us to investigate effectively the structure 
of the 
relative Brauer groups of binary cubic curves as we will see later.  
Notice that $\sqrt c\not\in k$, $\sqrt d\in k$ is equivalent to 
$c \in -3\,k^{*2}$, 
$\w\not\in k$, since $d=-27c$, and that $\sqrt c, \, \sqrt d\in k$ 
is equivalent to $\sqrt c, \, \w \in k$.

\proclaim{Corollary 5.3}  
	Let $k$ be a global field.  Let $\E\:Y^2=X^3+c$ and 
$\E'\: Y^2=X^3+d$ be the 
elliptic curves over $k$, where $d=-27c\in k^*$, and let 
$\a$, $\a'$ be the maps as above.	
\roster
\item "{\rm (i)}"  If $\sqrt c\not\in k$ but 
$\sqrt d\in k$ $($so, $\w\not\in k)$, 
then   
$$
|\im(\a)|\dd|\im(\a')|  \, = \,  3^{\rk+1}. \tag 5.3
$$ 
	If $\sqrt c, \sqrt d, \w\not\in k$, then 
$$
|\im(\a)|\dd|\im(\a')|  \, =  \, 3^{\rk}. \tag 5.4
$$ 

\item "{\rm (ii)}"  If $\sqrt c$, $\sqrt d\in k$  $($so, $\w \in k)$, then 
$$
|\im(\a) \, |= \, 3^{(\rk+2)/2}. \tag 5.5
$$  
	If $\sqrt c,\sqrt d\notin k$ but $\w \in k$, then 
$$
|\im(\a)| \, = \, 3^{\rk/2}. 
$$  
\endroster \endproclaim

\demo\nofrills{\smc {Proof.}\usualspace} \ (i) This is obvious by 
counting  the 
orders of $\ker\,(\l')$ and $\ker\,(\l)$, and  
applying Proposition~5.2.  

\noindent (ii) We first show that if $\w\in k$, then 
$\im(\a)=\im(\a')$.   To see 
this, notice that the equation $s^2=r^3-27c$ can be rewritten as 
$(\frac{s}{\sqrt{-3}^3})^2=(\frac{r}{\sqrt{-3}^2})^3+c$ since 
$-27=\sqrt{-3}^6$.   
This observation provides an isomorphism $\E(k)\cong \E'(k)$, 
given by 
$(r,s) \mapsto (\sqrt{-3}^2r,\sqrt{-3}^3s)$.  Moreover, 
we have 
$$
\textstyle\a'\big(\sqrt{-3}^2r,\sqrt{-3}^3s\big)  \,   = \ 
\ol{\sqrt{-3}^3s+\sqrt{-27c\,}\,} \, = \, \ol{s+\sqrt c\,}
 \, = \, \a(r,s).
$$  
This shows that $\im(\a)=\im(\a')$.  
	Now, if $\sqrt c,\;\sqrt d\in k$, then  
$|\ker\,(\l')|\dd|\ker\,(\l)| = 9$.  
Thus, by Proposition 5.2, we have 
$|\fim(\a)|^2\big/9=3^\rk$, which leads 
to (5.5).   Similarly, if $\sqrt c,\;\sqrt d\notin k$, 
then $|\im(\a)|^2 =3^\rk$. 
This completes the proof.    \qed\enddemo

\proclaim{Corollary 5.4}  
	Let $k$ be a global field and let $\E\:Y^2=X^3+c$ be the elliptic 
curve over $k$.  If  $\w \in k$, then $\E(k)$ has even rank. 
\endproclaim
	
	For the following proposition, let $k$ be a field with discrete 
valuation.  
For each prime spot $\frp$ of $k$, we  denote by $v_{\frp}$ the 
normalized discrete valuation corresponding to $\frp$.

\proclaim{Proposition 5.5}  
	Let $k$ be a field with discrete valuation.  Assume that $(r,s)$ 
is a $k$-rational point of the elliptic curve 
$\E\: Y^2=X^3+n^2,\ n\in k^*$.  
Then, for each prime spot $\frp$ of $k$ with 
$\vp(s+n) \not\equiv 0 \ (\mod\;3)$, 
one has $\vp(2n)\not\equiv 0 \ (\mod\;3)$.
\endproclaim

\demo\nofrills{\smc {Proof.}\usualspace}  
Since $r^3=s^2-n^2=(s+n)(s-n)$, it follows that 
	$$
\vp(s+n)+\vp(s-n) \, \equiv \,  0 \ (\mod\;3). \tag 5.6
$$  
This implies that $\vp(s+n) \ne \vp(s-n)$ since 
$\vp(s+n) \not\equiv 0 \ (\mod\;3)$.  
Thus, we have\break $\vp(2n)=\min\big(\vp(s+n),\vp(s-n)\big)$.  
Suppose now that 
$\vp(2n)\equiv 0 \ (\mod\;3)$.   Then 3 divides one of $\vp(s+n)$ 
and $\vp(s-n)$ 
but in fact both because of (5.6).  This contradicts the assumption 
that 
$\vp(s+n) \not\equiv 0 \ (\mod\;3)$.        \qed \enddemo

Suppose $k$ is the quotient field of a unique factorization domain
$R$.  For $a\in k^*$, we say that $b\in k^*$ is 
\lq\lq the" third-power-free part of $a$ if $b$ is third-power-free
in its prime factorization and $b \equiv a\ (\text{mod }k^{*3})$.
So, $b$ is unique up to $R^* \cap k^{*3}$. 
Proposition 5.5  shows that $\im(\a)$~lies in the 
$\zz/3\zz$-vector subspace of $k^*/k^{*3}$ generated 
by cube classes of units and all prime elements dividing the 
numerator or denominator of the third-power-free part of $2n$. 
In particular, we are most interested in 
the fields $\qq$ and $\qq(\w)$, which are the quotient fields of 
unique factorization 
domains $\zz$ and $\zz[\w]$ with units $\{\,\pm1\,\}$ and 
$\{\,\pm1,\pm\w,\pm\w^2\,\}$, 
respectively.  Since $-1\in {\qq}^{*3}\subseteq\qq(\w)^{*3}$, we 
have 

\proclaim{Corollary 5.6}  
	Let $\E\: Y^2=X^3+n^2$, $n\in k^*$ be an elliptic curve over a 
field $k$.  
Then, $\im(\a)$ lies in the $\zz/3\zz$-vector space 
generated by the following$:$
\roster
\item "{\rm (i)}"   all primes dividing the numerator or denominator 
of the third-power-free part of $2n$, if $k=\qq$.
\item "{\rm (ii)}"  $\w$ and all primes dividing the numerator or 
denominator of the third-power-free part of $2n$, if $k=\qq(\w)$.
\endroster 
\endproclaim

	Consider the maps $\a\:k\to k(\sqrt c)^*/k(\sqrt c)^{*3}$ and 
$\a'\:k\to k(\sqrt d)^*/k(\sqrt d)^{*3}$.  
For $x=c$ or $d$, we note that the canonical map 
$$
i_x\:k(\sqrt x)^*/k(\sqrt x)^{*3} \to 
k(\sqrt c,\sqrt  d)^*/k(\sqrt c,\sqrt  d)^{*3}
$$ 
is injective.  Thus, we will identify a cube class in 
$k(\sqrt x)^*/k(\sqrt x)^{*3}$ 
with its image under the map~$i_x$ and write $\im(\a')$ for 
$\im(i_x\!\circ\!\a')$ 
in $k(\sqrt c,\sqrt  d)^*/k(\sqrt c,\sqrt  d)^{*3}$ by abuse 
of notation.

\proclaim{Proposition 5.7}  
	Let $k$ be a field which does not contain $\w$.  Let $\E$ and 
$\E'$ be the 
elliptic curves over $k$ given by $\E\:Y^2=X^3+c$ and 
$\E'\: Y^2=X^3+d$, 
where  $d=-27c \in k^*$.  Then 
$$
\im(\a)\cap\im(\a') \, = \, \{\,\ol{1}\,\}.$$
\endproclaim

\demo\nofrills{\smc {Proof.}\usualspace}  
Let $\ol{\xi}\in\im(\a)\cap\im(\a')$.   Since $\ol{\xi}\in\im(\a)$, 
$\xi$ is of the form $s+\sqrt{c\,}$ for some $(r,s)$ in $\E(k)$. 
 Let 
$G=\Gal(k(\sqrt c,\sqrt d)/k)$ and take the $\s\in G$ with 
$\sqrt c \mapsto -\sqrt c$ but $\sqrt d$ fixed.  Then, $G$ acts on 
$k(\sqrt c,\sqrt d)^*/k(\sqrt c,\sqrt d)^{*3}$ and thus   
$$
\s(\ol{\xi}) \, = \, \ol{\s(s+\sqrt c)} \,  \, = \, 
 \, \ol{s-\sqrt c\,} \, \,  = \,  \, \ol{s+\sqrt c\,}^{\,\,-1}  
= \, \ol{\xi}^{\,\,-1},
$$ 
since $(s+\sqrt c)(s-\sqrt c)=r^3$.  On the other hand, since 
$\ol{\xi}\in\im(\a')$, 
it follows that $\ol{\xi}$ is fixed under the map $\s$.  
In other words, 
$\ol{\xi}=\ol{\xi}^{\,\,-1}$, which implies that $\ol{\xi}=\ol{1}$.     
\qed \enddemo

	Let $M$ be a finitely generated abelian group.  Recall that the 
rank of $M$ is defined by 
${\rank (M)=\dim_{\qq}(M\otimes_{\zz}\qq)}$. Assume that 
$G$ is a group of order 2 and that $G$ acts on $M$ 
(so~$M$ has a $\zz[G]$-module 
structure).   Denote by $M^G$ the submodule of $M$  fixed under 
the action of 
$G$ and by~$\widetilde M^G$  the submodule of $M$  fixed under the 
twisted 
action of $G$, that is, for the $\sigma\in G$, $\sigma\ne\id$, 
$$
\widetilde M^G  \, = \,  \{\,m\in M\,|\,\sigma(m)=-m\,\}.
$$ 

\proclaim{Lemma 5.8}  
	Let $M$ be a finitely generated abelian group.  With $M^G$ and 
$\widetilde M^G$ as above, we have 
$$
\rank(M) \, = \, \rank(M^G)+\rank(\widetilde M^G).
$$
\endproclaim

\demo\nofrills{\smc {Proof.}\usualspace}  
Define 
$$
M_1 \, = \, \{\,m+\sigma(m)\,|\,m\in M\,\}\  \ \ \text{and} 
\ \ \ M_2 \, = \, \{\,m-\sigma(m)\,|\,m\in M\,\}.
$$   
It is easy to check that $M_1$ and $M_2$ are $G$-fixed  submodules 
of $M$ and that $M_1\subseteq M^G$ and $M_2\subseteq \widetilde M^G$.   Since $2m=m+\sigma(m)+m-\sigma(m)$, it follows that $2M\subseteq M_1 + M_2 \subseteq M^G + \widetilde M^G\subseteq  M$.   Since $2M$ has the same rank as $M$, it is obvious that $M^G \oplus \widetilde M^G$ has the same rank as $M$.  Now, note that $M^G \cap \widetilde M^G$ is the set of 2-torsion elements.  Thus, $M^G \cap \widetilde M^G$  must be finite (as $M$ is finitely generated) and $\rank(M^G \cap \widetilde M^G)=0$.  From the fact that $$\rank(M^G+\widetilde M^G)+\rank(M^G \cap \widetilde M^G)=\rank(M^G)+\rank(\widetilde M^G),$$ we obtain  $$\rank(M)=\rank(M^G+\widetilde M^G) =\rank(M^G)+\rank(\widetilde M^G).   \qed $$   \enddemo

	The following proposition plays an important role to explicitly 
describe the relative 
Brauer groups of binary cubic curves over $\qq(\w)$.  We give a 
short and direct proof of 
it though this result is already known 
(see [Gr, Lemma 16.1.1, p.~49]).

\proclaim{Proposition 5.9} 
	Let $k$ be a global field which does not contain $\w$ and let 
$L=k(\w)=k(\sqrt {-3})$.  
Let $\E$ be the elliptic curve over $k$ given by $\E\:Y^2=X^3+c$. 
    Then, 
$$
\rank(\E(L)) \, = \, 2\,\rank(\E(k)).
$$
\endproclaim

\demo\nofrills{\smc {Proof.}\usualspace}  
	Let $G=\Gal(L/k)$.  It is obvious that $\E(L)^G=\E(k)$.  We claim 
that $\widetilde{\E(L)}^G \cong \E'(k)$.  For this, take the 
$\s\in G$, $\s\ne\id$.  
Since $\s(r,s)=-(r,s)$ for any $(r,s)\in \widetilde{\E(L)}^G$, 
it follows that 
$(\s(r),\s(s))=(r,-s)$.  This shows that $r\in k$ and 
$s\in\sqrt{-3}\,k$ and so 
we can find $r',\,s'\in k$ such that $r=r'\big/\sqrt{-3}^2$ and 
$s=s'\big/\sqrt{-3}^3$.  	
	Then, the condition $s^2=r^3+c$ is equivalent to  
 $(s')^2=(r')^3-27c$.  
	Define a map 
$$
\textstyle \varphi\:  \E'(k) \to  
\big\{\,\big(\frac{r'}{\sqrt{-3}^2},
\frac{s'}{\sqrt{-3}^3}\big)\in \E(L)\,\big|
\,r',s'\in k\,\big\}\cup\{\,\O\,\}
$$  
given by $\varphi\big((r',s')\big) = \big(\frac{r'}{\sqrt{-3}^2},
\frac{s'}
{\sqrt{-3}^3}\big)$ and $\varphi(\O) = \O$.  Obviously, the map 
$\varphi$ is  
bijective and the codomain of $\varphi$ is $\widetilde{\E(L)}^G$.  
Further, using the 
formula of the group operation $\oplus$ on $\E$ (cf.~\S2), a 
straightforward 
calculation shows that the map $\varphi$ is a group homomorphism.  
Hence, we have 
$\widetilde {\E(L)}^G \cong \E'(k)$ as claimed.  
	Now, since there exist only finitely many 2-torsion points, it 
follows from Lemma 5.8 
that 	
	$$
\rank(\E(L)) \, = \, \rank(\E(L)^G)+\rank(\widetilde{\E(L)}^G)
 \, = \, \rank(\E(k))+\rank(\E'(k)) \, = \, 2\,\rank(\E(k)).
$$   
Note here that $\rank(\E(k))=\rank(\E'(k))$ since the curves 
$\E$ and $\E'$ are isogenous.   This completes the proof.   
\qed \enddemo

	Proposition 5.9 confirms that the ranks of the elliptic curves in 
our situation are 
always even as we have seen in Corollary 5.4.     

\proclaim{Corollary 5.10} 
	Let $k$ be a global field not containing $\w$, and let 
$L=k(\w)=k(\sqrt {-3})$.  Let 
$\E$ be the elliptic curve over $k$ given by $\E\:Y^2=X^3+c$.     
Then, as $\zz/3\zz$-vector spaces, 
$$
\dim\,\a(\E(L)) \, = \, \dim\,\a(\E(k))+\dim\,\a'(\E'(k)).
$$
\endproclaim

	Let $\E$ be the elliptic curve defined by $\E\:Y^2=X^3+c$ over 
$\qq$, where $c$ is a sixth-power-free integer.  If $\E_{\tors}(\qq)$ denotes the torsion subgroup of 
$\E(\qq)$, then it is known (cf.\;[Hu, Th.~ 3.3, p.\;35]) that
        $$
\E_{\tors}(\qq)  \, \cong \ \left\{\,\aligned &\zz/6\zz  \ \ \  
\text{if} \ c=1,   \\
&\zz/3\zz  \ \ \ \text{if} \; c\ne 1 \text{ is a square, or } 
c =-432=-2^4\dd3^3,  \\
&\zz/2\zz  \ \ \ \text{if} \; c\ne 1 \text{ is a cube},  \\
&  \{\, \O\,\}  \hskip .17 in \text{otherwise.} \endaligned 
\right. \tag 5.7
$$

	In the curve $\E$ above, we may assume that $c$ is an integer 
because otherwise we can 
take an isomorphic curve  by multiplying $c$ by the sixth power of 
the denominator of $c$.  
We can now handle the rank 0 case completely.

\proclaim{Proposition 5.11}  
	 For $a,b\in {\qq}^*$, let $\C\: Z^3=aX^3+bY^3$ be a 
projective curve over 
$\qq$, and let $\E\: Y^2=X^3+c$ be the Jacobian of \; $\C$,  
where $c=-\frac{27}4a^2b^2$ is a 
sixth-power-free integer.  Let $k=\qq(\w)$.
\roster
\item "{\rm (i)}" If $c=-432$, then $\E(\qq)$ has rank $0$ and  
$$
\Br(\qq(\C)/\qq) \, = \,  \langle\,[(\qq(u)/\qq,\tau,a)]\,\rangle
\  \ \ \text{and}
 \ \ \ \Br(k(\C)/k)  \, = \,  \langle\,[(a,\w)_{\w}]\,\rangle,
$$ 
where $u \, = \, 
\cbrt{36+12\sqrt{-3\,}\,}+\,\cbrt{36-12\sqrt{-3\,}\,}$, with the 
cube roots chosen so that their product is $12$.  
\medskip
\item "{\rm (ii)}"  If $\E(\qq)$ has rank $0$ and $c\ne -432$, then 
$$
\Br(\qq(\C)/\qq) \,=\, \{\, 0\,\} \ \ \text{while} \ \ 
\Br(k(\C)/k) \,=\, 
\langle\,[(a,b)_{\w}]\,\rangle.
$$      
\endroster   \endproclaim

\demo\nofrills{\smc {Proof.}\usualspace} \ Since $\E(\qq)$ is 
assumed to have rank 0,  
we have $\E(\qq)=\E_{\tors}(\qq)$.  According to (5.7), all possible 
cases of rank 0 with $c=-\frac{27}4a^2b^2$ can be represented as 
follows:
        $$
\E(\qq)  \, \cong  \left\{\,\aligned 
&\zz/3\zz  \ \ \ \text{if} \; c=-432=-2^4\dd3^3,  \\
&\zz/2\zz  \ \ \ \text{if} \; c=-27,  \\
&  \{\, \O\,\}  \hskip .17 in \text{otherwise.} \endaligned \right. 
$$
	
	(i) \; If $c=-432$, then $\E(\qq)$ indeed has rank 0 (cf.\;[SAGE] 
or see Remark 5.12 
below) with rational points $(12,\pm 36)$ other than $\O$.    
A short calculation shows that 
$$
\a(12,36)\, =\, \ol{36+\sqrt{-3^3\dd4^2}\,}\, =\, 
\ol{4\dd3\sqrt{-3}(-\sqrt{-3}+1)}\, =\, 
\ol{-2^3\dd\sqrt{-3}^3(\textstyle\frac{-1+\sqrt{-3}}2)}\, 
=\, \ol{\w}.
$$  
Since the point $(12,-36)$ is the inverse of $(12,36)$, we 
obviously have 
$\a(12,-36)=\ol{\w^2}$, and thus $\a(\E(\qq))=\langle \w \rangle$.  
Further, 
since $\E(k)$ also has rank 0 by Proposition 5.9, it follows 
from (5.5) that 
$\a(\E(k))=\langle \w \rangle$.  Hence, by Propositions 4.1 and 
4.4, we have  
$$
\Br(\qq(\C)/\qq) \, =\,  \langle\,[(\qq(u)/\qq,\tau,a)]\,\rangle
\  \ \ \text{and}
\  \ \ \Br(k(\C)/k) \, =\,  \langle\,[(a,\w)_{\w}]\,\rangle,
$$ 
where $u\, =\, \cbrt{36+12\sqrt{-3\,}\,}+\cbrt{36-12\sqrt{-3\,}\,}$
for cube roots with product $12$, 
by Remark 4.5.

	(ii) \; If $c=-27$, then $\E(\qq)$ has rank 0 (cf.\;[SAGE]).  
Since $\E(\qq)
\cong \zz/2\zz$ but $\qq(\sc)^*/\qq(\sc)^{*3}$ is 3-torsion, we 
have $\a(\E(\qq))=\{\,\ol 1\,\}$.   If $\E(\qq)$ has rank 0 with 
$c\ne-432$ or $-27$, 
then $\E(\qq)=\{\,\O\,\}$ by (5.7) and thus 
$\a(\E(\qq))=\{\,\ol 1\,\}$.  On the other 
hand, since $\E(k)$ also has rank 0 and contains two $k$-rational 
points $(0,\pm\sc)$ 
other than $\O$, it follows that 
$\a(\E(k))=\langle \ol{2\sc}\rangle$.  Hence,  by 
Corollary 4.6 and Remark 4.2\,(i) respectively, we have   
$$
\Br(\qq(\C)/\qq)\, =\, \{0\}\ \ \ \text{and}\ \ \ \Br(k(\C)/k) 
\, = \, \langle\,[(a,b)_{\w}]\,\rangle. \ \ \ \qed
$$

\noindent{\bf Remark 5.12}.  \ In Proposition 5.11\,(i), we can 
 show directly that the 
rank of $\E(\qq)$ is 0.  To show this, for the curve 
$\E\:Y^2=X^3-3^3\dd4^2$, consider 
the isogenous curve $\E'\:Y^2=X^3+4^2$.  The image of the map 
$\a'\:\E'(\qq) \to 
\qq^*/\qq^{*3}$ must be trivial by Corollary 5.6\;(i).  On the 
other hand, in 
$k=\qq(\sqrt{-3})$, notice that $-3^3\dd4^2$ is a square and 
$\ol{\,2\dd\sqrt{-3^3\dd4^2\,}\,}=\ol{1}$ in $k^*/k^{*3}$.  By 
Corollary 5.6\;(ii), 
we then have $\a(\E(k))\subseteq \langle\,\ol{\w}\,\rangle$.  
Since 
$\E(\qq)\subseteq \E(k)$, it follows that 
$\a(\E(\qq))\subseteq \langle\,\ol{\w}\,\rangle$.   
Thanks to the formula in (5.3), the rank of $\E(\qq)$ must be 0.

\bigskip                                         


	For the rest of this section, we provide various examples of 
relative Brauer groups of 
function fields of diagonal cubic forms.  For a curve 
$\E\:Y^2=X^3+c$, we may take as its 
isogenous curve $\E'\:Y^2=X^3-\frac{\textstyle c}{27}$ instead of 
$Y^2=X^3-27c$ since 
these two curves are isomorphic.   \smallskip

\noindent{\bf Example 5.13}.  \ Consider the projective curve 
 $\C\: Z^3=aX^3+bY^3 \;(a,b\in\qq^*)$ with $ab=2\dd11$.  
The Jacobian of $\C$ is $\E\:Y^2=X^3-27\dd11^2$ and its 
isogenous curve is 
$\E'\:Y^2=X^3+11^2$.  It is known (cf.\; [SAGE]) that the curve $\E'$ (so $\E$ as well) 
has rank 1. 
	Notice that $\E'$ has integer points $(0,11)$, $(12,43)$, which are 
sufficient to 
determine $\a(\E(\qq))$ and $\a'(\E'(\qq))$.  The images of these 
points in 
$\qq^*/\qq^{*3}$ are $\ol{2\dd11}$, $\ol {2\dd27}\,(= \overline{2})$.  
	Since $\a'(\E'(\qq))$ is contained in the $\zz/3\zz$-vector space 
generated by all 
primes  dividing the third-power-free part of $2n$ by 
Corollary 5.6\,(i), it follows 
that $\a'(\E'(\qq))=\big\langle\, \ol{2}, \ol{11} \,\big\rangle$.  
Then, the formula in 
(5.3) with the rank information implies that 
$\a(\E(\qq))=\{\,\ol{1}\,\}$.  Hence, by 
Corollary 4.6, we have  
$$
\Br(\qq(\C)/\qq) \,= \, \{\,0\,\}.
$$

	For $k=\qq(\w)$, we next consider the  curves $\E$ and $\E'$ over 
$k$.   Note that the 
primes 2 and 11 are unramified in $k$.  Since the two integer 
points of the curve $\E'$ 
above again have images $\ol{2\dd11}$, $\ol {2}$ in $k^*/k^{*3}$, 
this implies that  
$\langle\,\ol {2},\ol {11}\,\rangle \subseteq\a(\E(k))\subseteq
\langle\,\ol{\w},\,\ol {2},\ol {11}\,\rangle $. 
	However, since $\E(k)$ has rank 2 by Proposition 5.9, it follows 
that 
$\a(\E(k))=\langle\,\ol {2},\,\ol {11}\,\rangle$ and  
$$
\Br(k(\C)/k)\, =\, \langle\,[(a,2)_{\w}],\,[(a,11)_{\w}]\,\,\rangle.
$$
This is valid for any choice of $a\in \qq^*$ when we set $b = 22/a$.
Clearly $\big|\Br(k(\C)/k)\big|$ depends on the choice of $a$.

	For a specific example, let $a=2$ (and so $b=11$).  Note that 
$(2,2)_{\w}\cong(2,-1)_{\w}$ is obviously split and 
that $(2,11)_{\w}$ is also split 
since 11 is a norm for the extension $k(\cbrt{2})/k$, as 
$11=3^3-2\dd2^3$.  
Hence, we have 
$$
\Br(\qq(\C)/\qq)\,  =\, \Br(k(\C)/k) \, =\, \{\,0\,\}.
$$  

\medskip

\noindent{\bf Example 5.14}.  \   \ 
Consider the projective curve  $\C\: Z^3=aX^3+bY^3$ with 
$ab=2\dd3\dd5$.  The 
Jacobian of $\C$ is $\E\:Y^2=X^3-27\dd15^2$ and its isogenous 
curve is 
$\E'\:Y^2=X^3+15^2$.  It is known (cf.\; [SAGE]) that the curve 
$\E'$ has rank 2. 
 Notice that $\E'$ has integer points $(-5,-10)$, $(0,15)$, 
$(4,-17)$, whose images 
in $\qq^*/\qq^{*3}$ are $\overline {5}$, 
$\overline {30}\,(= \overline{2\!\cdot\!3\!\cdot\!5}$), 
$\overline {-2}\,(=\overline {2})$, respectively.  
Then, Corollary~5.6\,(i) tells us that  
$\a'(\E'(\qq))= \langle\,\overline {2},\,
\overline {3},\,\overline {5}\,\rangle$.  Thanks to the formula 
in (5.3), we have 
$\a(\E(\qq))=\{\,\ol1\,\}$.  Hence, by Corollary 4.6,   
$$
\Br(\qq(\C)/\qq)\,  =\,  \{\, 0\,\}.
$$  

	Next, we observe that the primes 2 and 5 are unramified and 
3 is ramified in $k=\qq(\sqrt{-3})$ 
(note that $\sqrt{-3}$ is the prime of $k$ above 3).  The 
three integer  points of the 
curve $\E'$ above again have images  $\ol {5}$,  
$\ol {30}$, $\ol {2}$ in $k^*/k^{*3}$ 
respectively.  Since   $30=2\dd\sqrt{-3}^2\dd 5$, it follows from 
Corollary 5.6\,(ii) 
that  
$$
\langle\,\ol {2},\,\ol {\sqrt{-3}\,},\,\ol {5}\,\rangle 
\, \subseteq\, \a(\E(k))\, =\, \a'(\E'(k))\, 
\subseteq\, \langle\,\ol{\w},\,\ol {2},\,
\ol {\sqrt{-3}\,},\,\ol {5}\,\rangle .
$$  
	Since $\E(\qq)$ has rank 2, it follows from Proposition 5.9 
that $\E(k)$ has rank 4.  
Hence, we have 
$\a(\E(k))=\langle\,\ol {2},\,\ol {\sqrt{-3}\,},\,\ol {5}\,\rangle$ 
and  
$$
\Br(k(\C)/k)\, =\, 
\langle\,[(a,2)_{\w}],\,[(a,\sqrt{-3})_{\w}],\,[(a,5)_{\w}]\,\,
\rangle.
$$

\medskip

\noindent{\bf Example 5.15}.  \ Take for  $\C\:Z^3=aX^3+bY^3$ with 
$ab=2\dd13$.  The 
Jacobian of $\C$ is\break $\E\:Y^2=X^3-27\dd13^2$ and its isogenous 
curve is 
$\E'\:Y^2=X^3+13^2$.  By [SAGE], we see that the curve $\E'$ has 
rank 1.  Further, 
it can be checked that $\E$ and $\E'$ have integer points 
$(39,234)$ and $(0,13)$, 
respectively.  Note that 13 is split in $k=\qq(\w)$; in fact, 
$13=\frp\frq$ where 
$\frp=1+2\sqrt{-3}$ and $\frq=1-2\sqrt{-3}$.  
	We now claim that $\a(\E(\qq))=\langle\,\ol{\frp\frq^2}\,\rangle$.  
To see this, 
observe that 
$$
234+\sqrt{-27\dd13^2}\, =\, 3\dd13\dd\sqrt{-3}(-2\sqrt{-3}+1)\, 
=\, 
-\sqrt{-3}^3\dd13\dd(1-2\sqrt{-3})\, =\, -\sqrt{-3}^3\dd\frp\frq^2.
$$  
	Since $-\sqrt{-3}^3$ is a cube in $k$, it follows that 
$\langle\,\ol{\frp\frq^2}\,\rangle \subseteq \a(\E(\qq))$.  On the 
other hand, 
we know $\langle\,\ol{2\dd13}\,\rangle \subseteq \a'(\E'(\qq))$.   
However, since 
the rank of the curve $\E$ is 1, it follows from the formula 
in~(5.3) that 
$$
\a(\E(\qq))\, =\, \langle\,\ol{\frp\frq^2}\,\rangle\, \cong\, 
 \zz/3\zz \text{\  \ \ and\  \ \ } 
\a'(\E'(\qq))\, =\, \langle\,\ol{2\dd13}\,\rangle\, \cong\,  \zz/3\zz.
$$
	Hence, by Proposition 4.4, 
$$
\Br(\qq(\C)/\qq)\, =\, \langle\,[(\qq(u)/\qq,\tau,a)]\,\rangle
$$ 
where 
$u\, =\, \cbrt{234+39\sqrt{-3\,}\,}\, +\,\, 
\cbrt{234-39\sqrt{-3\,}\,}$\!, with cube
roots chosen so that their product is~$39$. 
 Note here that $(\qq(u)/\qq,\tau,a)\otimes_{\qq} k 
\cong (a,\frp\frq^2)_{\w}$.

	Next, over $k$, we have 
$\langle\,\ol{2\frp^2},\ol{2\frq^2}\,\,\rangle\subseteq 
\a(\E(k))\subseteq 
\langle \,\ol{\w},\,\ol{2},\,\ol{\frp},\,\ol{\frq}\,\,\rangle$.   
However, since $\E(k)$ has rank 2 by Proposition 5.9, we have 
$\a(\E(k))=\langle\,\ol{2\frp^2},\ol{2\frq^2}\,\,\rangle$ and thus 
$$
\Br(k(\C)/k)\, =\, \langle\,[(a,2\frp^2)_{\w}],
\,[(a,2\frq^2)_{\w}]\,\,\rangle.
$$

	For a specific example, let $a=3$, so we have \; 
$\C\:Z^3=3X^3+\frac{2\cdot13}3Y^3$.  
In what follows, we describe the nonsplit algebras in the relative 
Brauer groups by 
means of local invariants.  For this, it suffices to consider the 
symbol algebras 
over $k$ locally at the prime spots 2,~$\sm$,~$\frp$,~and~$\frq$, 
where 
$\frp\dd\frq=13$.  First, over the dyadic completion $\k2$, the 
symbol algebra 
$(3,2\frp^2;\k2)_\w \cong(3,2;\k2)_\w$ is split since 3 is a cube  
modulo 2 and by 
Hensel's lemma.  
	Next, it is obvious that the invariant of $(3,2\frp^2;\kq)_\w$ at 
$\frq$, 
denoted $\inv_\frq(3,2\frp^2;\kq)_\w$, is trivial in $\qq/\zz$ 
since 2, 3, and 
$\frp$ are units $\frq\text{-adically}$.  
	In order to determine $\inv_\frp(3,2\frp^2;\kp)_\w$, consider the 
algebra 
$(3,\frp;\kp)_\w$.  Since 13 is split in $k$, notice that 
$\kp=\q{13}$ and the 
residue field $\ol k_{\frp}=\zz_{13}$.  Inside the algebra, we 
have defining 
relations \; $i^3=3,\; j^3=\frp,\;ij=\w ji$.  From the unramified 
extension 
$\kp(i)/\kp$ of degree~3, we consider the induced automorphism 
$\varphi$ of 
the residue field $\zz_{13}(\ol \imath)$, which is some power of 
the Frobenius automorphism.  
	Since $\frp=1+2\sm$, it follows that 
$\sm\equiv \frac{-1}2\;(\mod\;\frp)$.  
Observing that 
${\ol\w=\ol{(-1+\sm)\,/\,2\,}=\ol{(-1-(1/2))\,/\,2\,}=
\ol{-3/4}=
\ol{\,3^2}}$  in $\zz_{13}$ and so $\ol{\w^{-1}}=\ol3$, we have 
$\varphi(\ol i)=\ol{jij^{-1}}= \ol{\w^{-1}i} =\ol{3i}=
\ol i^{\,13}$.
	Therefore the map $\varphi$, conjugation by $\ol j$, is 
the Frobenius automorphism 
and so $\inv_\frp(3,2\frp^2;\kp)_\w = \ol{\,2/3\,}$.  By 
Hilbert's Reciprocity Law, 
the symbol algebra $(3,2\frp^2;k)_\w$ has local invariant 
$\frac23$ at $\frp$, 
$\frac13$ at $\sm$, and 0 at all other prime spots.
	By similar calculations, it can be shown that the algebra 
$(3,2\frq^2;k)_\w$ has 
local invariant $\frac13$ at $\frq$, $\frac23$ at $\sm$, and 0 
at all other prime spots.
	As $(3,\frp\frq^2)_\w \sim (3,2\frp^2)_\w^2\ok (3,2\frq^2)_\w$, 
the algebra 
$(3,\frp\frq^2)_\w$ has local invariant $\frac13$ at each of $\frp$, 
$\frq$, and $\sm$.
	Consequently, the cyclic algebra $A:=(\qq(u)/\qq,\tau,2)$ satisfying 
$A\otimes_{\qq} k \cong (a,\frp\frq^2)_{\w}$ has local 
invariant $\frac13$ 
at 13 and $\frac23$ at 3, and 0 at other prime spots.

\medskip

\noindent{\bf Example 5.16}.  \ Consider $\C\: Z^3=aX^3+bY^3$ with 
$ab=136=2^3\dd17$.  
The  Jacobian of $\C$ is $\E\:Y^2=X^3-27\dd(4\dd17)^2$ and its
 isogenous curve 
is $\E'\:Y^2=X^3+(4\dd17)^2$.  In this example, we can easily show 
that the 
rank of $\E(k)$ is 2 (so, the rank of $\E(\qq)$ is 1).  Notice 
that $\E$ and $\E'$ 
have integer points $(84,684)$ and $(0,68)$, respectively, and 
that 17 is 
unramified in $k$.  Since $2\dd68=2^3\dd17$, it follows  from 
Corollary 5.6\;(i) 
that $\a'(\E'(\qq))=\langle\,\ol{17}\,\rangle$.  
	On the other hand, observe that 
$$
\textstyle 684+\sqrt{-27\dd4^2\dd17^2}\, \, =\,\,  
3\sm\cdot4(-19\sqrt{-3}+17)\, \, =\, \, 
\sqrt{-3}^3\cdot2^3\cdot\frac{19\sm-17}{2} 
\, \, =\, \,  -\sqrt{-3}^3\cdot2^3\cdot(2-\sqrt{-3})^3\w^2.
$$ 
	In other words, $\ol{684+\sqrt{-27\dd4^2\dd17^2}\,}=\ol{\,\w^2}$, so 
$\langle\,\ol{\w}\,\rangle \subseteq \a(\E(\qq))$, which shows that 
$\a(\E(k))=\langle\,\ol{\w},\,\ol{17}\,\rangle $  
by Corollary 5.6\;(ii).  
Therefore, the rank of $\E(k)$ must be 2 by the formula in (5.5).   
Since 
$\ol{17}\not\in \E(\qq)$ by Proposition 5.7, we have 
$\a(\E(\qq))=\langle\,\ol{\w}\,\rangle$.   
	Hence, we have, 
$$
\Br(\qq(\C)/\qq)\, =\, \langle\,[(\qq(u)/\qq,\tau,a)]\,\rangle
$$ 
where $u\, =\, \cbrt{684+204\sqrt{-3\,}\,}+\, \,
\cbrt{684-204\sqrt{-3\,}\,}$, with cube
roots chosen so that their product is~$84$.  Note here 
that $(\qq(u)/\qq,\tau,a)\otimes_{\qq} k \cong 
(a,\w)_{\w}$.   Moreover, 
we have 
$$
\Br(k(\C)/k)\, =\, \langle\,[(a,\w)_{\w}],\,[(a,b)_{\w}]\,\rangle,
$$
since $(a,17)_{\w}\cong (a,b)_{\w}$.
 
\medskip

\noindent{\bf Example 5.17}.  \ Consider $\C\: Z^3=aX^3+bY^3$ with 
$ab=728=2^3 \dd 7 \dd 13$.  
The  Jacobian of $\C$ is $\E\:Y^2=X^3-27(4\dd7\dd13)^2$ and its 
isogenous curve 
is $\E'\:Y^2=X^3+(4\dd7\dd13)^2$.  It is known (cf.\;[SAGE]) that 
the curve $\E$ has 
rank 2 over $\qq$. Also, $\E$ has rational points $(156,468)$ and 
$(196,1988)$ (found by computer search) and  
$\E'$ has   rational point $(0,364)$. 
Note that both 7 and 13 are split in $k=\qq(\w)$; in fact, 
$7=\frp_1\frp_2$ where 
$\frp_1=2+\sqrt{-3}$ and $\frp_2=2-\sqrt{-3}$, and 
$13=\frq_1\frq_2$ where 
$\frq_1=1+2\sqrt{-3}$ and $\frq_2=1-2\sqrt{-3}$. 
	Since $2\dd364=2^3\dd7\dd13$, it follows that 
$\langle\,\ol{7\dd13}\,\rangle \subseteq \a'(\E'(\qq))$ by 
Corollary 5.6\;(i).  
On the other hand, observe that 
	$$
468+\sqrt{-27(\dd4\dd7\dd13)^2}  \,\,  = \, 
\, 4\dd9\dd13+4\dd3\dd7\dd13\sqrt{-3} 
 \,\,  =\, \, 2^3\dd\sqrt{-3}^3\dd13\dd\frac{\sqrt{-3}-7}{2} 
\, \, =\, \,  2^3\dd\sqrt{-3}^3 \dd 13\frq_2\w^2.
$$   
Similarly, it can be shown that 
$1988+\sqrt{-27(\dd4\dd7\dd13)^2} = 2^3\frp_1^5\frp_2\w^2$.  So,
${\langle\,\ol{\frq_1\frq_2^2\w^2},\,\ol{\frp_1^2\frp_2\w^2}\rangle 
\subseteq \a(\E(\qq))}$.  
	Since the rank of $\E(\qq)$ is 2, the formula in (5.3) shows that 
$\a(\E(\qq)) = \langle\,\ol{\frq_1\frq_2^2\w^2},\, 
\ol{\frp_1^2\frp_2\w^2}\rangle$.  
Likewise, since the rank of $\E(k)$ is 4, we have 
$\a(\E(k))=
\langle\,\ol{7\dd13},\, \ol{\frp_1^2\frp_2\w^2},
\, \ol{\frq_1\frq_2^2\w^2}\rangle$ by (5.5). 
 
	Hence, we conclude that  
$$
\Br(\qq(\C)/\qq)=\big\langle\,[(\qq(u_1)/\qq,\tau_1,a)],
[(\qq(u_2)/\qq,\tau_2,a)]\,\big\rangle
$$ 
where
$$
u_1\,=\,\,\cbrt{468+1092\sqrt{-3\,}\,}\,+\,\,
\cbrt{468-1092\sqrt{-3\,}\,} \ \ 
 \text{and} \  \ 
u_2\,=\,\,\cbrt{1988+1092\sqrt{-3\,}\,}\,+\,\,\cbrt{1988-1092\sqrt{-3\,}\,}.
$$
The cube roots in the formula for $u_1$ are chosen so that their
product is real; likewise for $u_2$.
 Also, the automorphisms 
$\tau_1$ and $\tau_2$ are chosen so that  
${(\qq(u_1)/\qq,\tau_1,a)\otimes_{\qq} k  \cong 
(a,\frp_1^2\frp_2\w^2)_{\w}}$ and
${(\qq(u_2)/\qq,\tau_2,a)\otimes_{\qq} k \cong 
(a,\frq_1\frq_2^2\w^2)_{\w}}$. 
        Furthermore, we have
$$
\Br(k(\C)/k)=\langle\,[(a,7\dd13)_{\w}],[(a,\frp_1^2\frp_2\w^2)_{\w}],\,[(a,\frq
_1\frq_2^2\w^2)_{\w}]\,\rangle.
$$

\bigskip\bigskip


\noindent {\bf \S6. \ The nondiagonal case}            \medskip


\def\tb{\ol\tau}
\def\tf{\widetilde f}
\def\tt{\widetilde\tau}
\def\ta{\widetilde a}
\def\tb{\widetilde b}

	Let $\C =\C_f$ be the smooth projective genus 1 curve given by 
 $$
\C_f\: Z^3\, =\, f(X,Y)  \ \ \  \text{where} \  \ 
f(X,Y)\, =\, AX^3+3BX^2Y+3CXY^2+DY^3,  
\tag 6.1
$$ 
with $A,B,C,D\in k$.  In previous sections, we have described the 
relative 
Brauer groups $\Br(k(\C)/k)$ in the diagonal case, that is, when
 $B=C=0$.  We now turn to the nondiagonal case. 

	If $f$ is diagonalizable by a linear change of variables over $k$, 
then we can apply the diagonal case to compute $\Br(k(\C)/k)$.  
We recall when such a diagonalization is possible, following the 
approach in [D,\;pp.\,16-17].

	The Hessian matrix of $f$ is 
$$
H_f(X,Y)  \ = \  \left(\matrix{\partial^2 f}/{\partial X^2} 
& {\partial^2 f}/{\partial X\,\partial Y} \\ 
{\partial^2 f}/{\partial Y\partial X} & 
{\partial^2 f}/{\partial Y^2}  \endmatrix\right) 
 \ = \  6\left(\matrix AX+BY & BX+CY \\ BX+CY  & CX+DY  
\endmatrix\right). 
$$
	The Hessian determinant of $f$ is defined to be 
$$
h_f(X,Y) \, = \,  \det\big(H_f(X,Y)\big)  \, =
 \  36(RX^2+2SXY+TY^2),
$$ 
where 
$$
R \, = \, AC-B^2, \quad 2S \, = \, AD-BC, \quad \text{and} 
\quad T \, = \, BD-C^2.  \tag 6.2
$$  
The discriminant $\D_f$ of $f$ is 
$\frac1{\,72^2}\dd\text{(discriminant of } h_f)$, that is, 
$$
\D_f \, = \, S^2-RT  \, = \,  \,  
\textstyle \frac14\big(A^2D^2-3B^2C^2+4AC^3+4B^3D-6ABCD \big).   
\tag 6.3
$$ 
	\big(Note that $\D_f$ is not the discriminant of the cubic 
polynomial $f(X,1)\!\in\! k[X]$.  Indeed,\break 
$\D_f= \textstyle \frac{-1}{\,4\cdot 27\,}\disc\,(f(X,1)) 
\equiv -\,3\,\disc\,(f(X,1))  \ (\mod\;k^{*2}).\big)$ 
	We assume throughout that $f$ is nondegenerate, that is, 
$\D_f\ne 0$.  This is equivalent to the condition that 
$f(X,1)$ has 3 distinct roots in an algebraic closure of $k$.

	A short computation shows that $f$ is diagonal (i.e., $B=C=0$) 
if and only if  $R=T=0$.  
In this case, 
$\D_f=({AD}/2)^2\in k^{*2}$.  
We can use this to see that $f$ is 
diagonalizable over $k$ if and only if the quadratic 
form $h_f(X,Y)$ is isotropic over $k$, that is, if
and only if its 
discriminant, $\Delta_f$, is a square in $k$:  
Given new variables $U$ and $V$, we set
$X(U,V) =\a U+\b V$ and $Y(U,V) =\c U+\d V$ with 
$\a, \b, \c, \d \in k$ chosen so that the matrix 
$Q=\left(\smallmatrix \a & \b \\ \c & \d \endsmallmatrix\right)$ 
satisfies $\det(Q) \ne 0$.  Let $\epsilon = \det(Q)$.  
Then, for  $\tf (U,V)=f(X(U,V),Y(U,V))$, the 
chain rule yields  
$$
H_{\wt{f}}(U,V) \ = \ Q^t\, H_f(X(U,V),Y(U,V))\, Q.
$$  
Hence, $h_{\wt{f}}(U,V)=\e^2h_f(X(U,V),Y(U,V))$.  
	If we write   
$$
h_f  \, = \,  36(RX^2+2SXY+TY^2) \ =
 \ 36(\ol R U^2+2\ol S UV+\ol T V^2),
$$ 
then 
$$
\left(\smallmatrix \ol R & \ol S \\
 \ol S & \ol T \endsmallmatrix\right)
\, =\, \, Q^t
\left(\smallmatrix R & S \\ S & T \endsmallmatrix\right)Q,
$$ 
and thus  
$\ol R\,\ol T - \ol S^2=\e^2(RT-S^2)$.  
	By expressing
$h_{\tf} = 36(\wt R U^2+2\wt S UV +\wt TV^2)$,  since\break 
$h_{\tf}(U,V)=\e^2h_{f}(X(U,V),Y(U,V))$, it follows that 
$\wt R=\e^2\ol R$,  $\wt S=\e^2\ol S$, and $\wt T=\e^2\ol T$.  
	Hence, 
$$
\D_{\tf} \, =  \, \wt S^2-\wt R\wt T \, = \, 
\e^4(\ol S^2-\ol R\,\ol T) \, = \, \e^6\D_f.
$$
	Thus, if $\D_f\not\in k^{*2}$, then $\D_{\tf}\not\in k^{*2}$ and so 
$\tf$ is not diagonal over $k$.  Conversely, suppose 
that 
$\D_f\in k^{*2}$ but $f$ is not diagonal (so $R\ne 0$ or $T\ne 0$).  
If $R\ne 0$, we have 
$$
\textstyle  h_f \, = \,  \, 
\frac{\,36\,}R\big(RX+(S+\sqrt{\botsmash{\D_f}}\,)Y\big)
\big(RX+(S-\sqrt{\botsmash{\D_f}}\,)Y\big).
$$  
	Set $U= RX+(S+\sqrt{\botsmash{\D_f}}\,)Y$ and 
$V= RX+(S-\sqrt{\botsmash{\D_f}}\,)Y$ 
\ \(so $Q=\frac1{2R\sqrt{\botsmash{\D_f}}}
\left(\smallmatrix \sqrt{\botsmash{\D_f}} \  - \, S & \ 
\sqrt{\botsmash{\D_f}} \ + \, S \\
 \ R &  \  \ -R \endsmallmatrix\right)$  and 
$\e=-1\big/(2R\sqrt{\botsmash{\D_f}})\, \).$
Thus, $X = \frac1{2R\sqrt{\botsmash{\D_f}}}\big(
(\sqrt{\botsmash{\D_f}} \  - \, S) \,U+  
(\sqrt{\botsmash{\D_f}} \ + \, S )\,V\big)$
and $Y = \frac1{2\sqrt{\botsmash{\D_f}}}(U-V)$. 
Let $\wt f(U,V)\, =\, f(X(U,V),Y(U,V))$ as above. 
Because $h_f(U,V)=\frac{36}R UV$ we see that 
 $\wt R=\e^2\ol R=0$ and 
$\wt T=\e^2\ol T=0$, 
showing that $\wt f$ is diagonal. 
Indeed $\wt f(U,V)= a\,U^3 +b\,V^3$, 
where 
$$
\textstyle a \,=\, \wt f(1,0) \,=\, 
f\big(\frac{\vphantom{\D_f}\sqrt{\botsmash{\D_f}} \, \,- \, S}
{2R\sqrt{\botsmash{\D_f}}}, \frac1{2\sqrt{\botsmash{\D_f}}}\big)
\ \  \  \ \text{and}\ \  \  \ b\, = \, \wt f(0,1) \ = \ 
f\big(\frac{\vphantom{\D_f}\sqrt{\botsmash{\D_f}} \  + \, S}
{2R\sqrt{\botsmash{\D_f}}}, \frac{-1}{2\sqrt{\botsmash{\D_f}}}\big). 
$$

The case 
where $R=0$ (so $T\ne 0$) is treated analogously by putting 
$U= (S+\sqrt{\botsmash{\D_f}}\,)X+TY $ and 
$V= (S-\sqrt{\botsmash{\D_f}}\,)X+TY.$

\medskip

	Now, as in earlier sections, for $c\in k^*$ and  $\D=-c/27$, let 
$\E$ be 
the projective elliptic curve with affine model $\E\:Y^2=X^3+c$, 
and let 
$\T$ be its subgroup $\{\,\O,(0,\pm\sc)\,\}\subseteq \E(k_s)$.

\proclaim{Theorem 6.1}  
	The curves in the image of the canonical map 
$\Psi\:H^1(k,\T)\ra H^1(k,\E)$ are all the curves $\C_f$ as in 
$(6.1)$ with 
$\D_f\equiv \D\;(\mod\;k^{*6})$.
\endproclaim
\demo\nofrills{\smc {Proof.}\usualspace} \  
(We thank D.~Krashen for suggesting the corestriction argument 
given at the end of the proof.)       
	Assume first that $\D \in k^{*2}$.  We saw in the proof of 
Proposition 2.1 that  
the curves $\C_t$ in $\im(\Psi)$ have the form 
$\C_t\:X^3-tY^3=-54\,t^2\sd\,Z^3$ as $t$ ranges over $k^*$. This 
$\C_t$ is 
the curve $\C_f$ for $f(X,Y)=-\frac1{\,54\,t^2\sd\,}(X^3-tY^3)$; 
therefore, 
$$\D_f  \,  \, =  \, \,   {(-t)^2}/[\,4(-54\,t^2\sd)^4\,]  
   \, =  \,  \, {\D}[(18\,t\sd)^{-6}]  \,  \, \equiv
 \,  \,  \D\;(\mod\;k^{*6}).
$$  
	Now, for any $f$ as in (6.1), suppose that 
$\D_f\equiv \D\;(\mod\;k^{*6})$.  
We then have $\D_f\in k^{*2}$ by hypothesis.  Thus, by changing 
variables 
we may diagonalize $f$, reducing to the case where $B=C=0$; as 
we saw above, 
this change of variables preserves $\D_f$ modulo $k^{*6}$.   
	By Corollary 2.2, every $\C_f$ with $f$~diagonal is isomorphic 
to some $\C_t$.  
This completes the proof when $\D\in k^{*2}$.

	Next,  assume that $\D \notin k^{*2}$.  Let $L=k(\sd)$ and let 
$\tau$ be 
the nonidentity $k$-automorphism of $L$.  As we saw in (2.6) above, 
$H^1(k,\T)\cong \{\,\ol t\in L^*/L^{*3}\,\big|\,\tau(\ol t)
={\ol t}^{\,-1}\,\}$.  Take any $t\in L^*$ with 
$\tau(t)\equiv t^{-1}\;(\mod\;L^{*3})$.  
	Then $t^2/\tau(t)^2=t\big[t^3/(t\,\tau(t))^2\big]\equiv t\;
(\mod\;L^{*3})$; 
so, by replacing~$t$ by $t^2/\tau(t)^2$, we may assume that 
$\tau(t)=t^{-1}$.  Let 
$[\g_t]\in H^1(k,T)$ be the cohomology class corresponding 
to~$tL^{*3}$.  
	Since $\D\in L^{*2}$, Proposition 2.1 applies over $L$ and 
tells us that 
the curve over $L$ associated to $\Psi_L(\res_{k\to L}[\g_t])$ is 
$\C_{t,L}$ with equation 
$$
\C_{t,L}\:U^3-tV^3 \, = \, -54\,t^2\sd\,Z^3.   \tag 6.4
$$
	To see that $\C_{t,L}$ is actually defined over $k$, we make 
the change of variables 
$$
\textstyle X\, =\, \frac12(U+t^{-1}V) \ \  \text{and} 
\ \  Y\, =\, \frac{\sd}2(U-t^{-1}V),
\ \ \text{so} \ \ 
\textstyle U\, =\, \frac1{\sd}(\sd\,X+Y) \  \ \text{and}
  \ \ V\, =\, \frac t{\sd}(\sd\,X-Y).$$
	Then equation (6.4) becomes 
$$
\textstyle -54\,t^2\sd\,Z^3 \ = \ U^3-tV^3 \ =
 \ \D^{\,-3/2}\big[(\sd\,X+Y)^3-t^4(\sd\,X-Y)^3\big],
$$  
which simplifies to:  
$$
Z^3  \, \,  = \  AX^3+3BX^2Y+3CXY^2+DY^3,  
$$   
where 
$$
\textstyle A \, = \, \frac{\;t^2-t^{-2}\,}{54\sd}, \  \ 
B \, = \, \frac{\,-(t^2+t^{-2})\,}{54\D}, \  \ 
C \, = \, \frac{\;t^2-t^{-2}\,}{54\D\sd}, \  \ \text{and} \  \ 
D \, = \, \frac{\,-(t^2+t^{-2})\,}{54\D^2}.
\tag 6.5
$$
	Since $\tau(\sd)=-\sd$ and $\tau(t)=t^{-1}$,  the constants 
$A, B, C, D$ 
are each fixed by $\tau$, so lie in~$k$.  
Let $g(X,Y)= AX^3+3BX^2Y+3CXY^2+DY^3$. 
Then,   the curve 
${\C_g\:Z^3=g(X,Y)}$ is defined over $k$,
and $\C_g\times_k L \cong \C_{t,L}$ 
over $L$.  By [An, (3.4), (3.5), (3.8)], the Jacobian of 
$\C_g$ is the elliptic curve $\E$, 
which implies that $\C_g$ is the curve determined by some~ 
${\th\in H^1(k,\E)}$.  
Hence, 
$$
\res_{k\to L}(\th)  \,  = \,  \Psi_L(\res_{k\to L}[\g_t]) \ = \ 
\res_{k\to L}(\Psi[\g_t]) \ \ \text{in} \ H^1(L,\E).
$$
	Arguing as in the proof of Proposition 2.1, since $\th$ has 
points in a degree 3 
field extension of $k$, it follows that $\th$ is 3-torsion
in $H^1(k,\E)$.  
Also $\Psi[\gamma_t]$ is $3$-torsion, as $|\T| = 3.$
Because 
$[L\:\!k]=2$,  the restriction map $H^1(k,\E)\to H^1(L,\E)$ is 
injective 
on 3-torsion.  Hence, we have $\th=\Psi[\g_t]$, showing that 
the curve 
$\C_g$ lies in $\im(\Psi)$.  Furthermore,  as 
$(t^2-\,t^{-2})^2-(t^2+\,t^{-2})^2= -4$,  for the $R, S, T$ 
of (6.2) for this $g$ we have 
$$
 R   \, \,= \,  -(27\D)^{-2}, \quad  S \, = \, 0,\ \ \text{and} 
\quad
T \, =  \,  {27^{-2} \D^{-3}},
$$  
 so $\D_g=-RT=(9\D)^{-6}\D\equiv \D \; (\mod\; k^{*6})$.
	Since the choice of $[\g_t]\in H^1(k,T)$ was arbitrary, we have 
shown that 
every curve in $\im(\Psi)$ has the specified form.

	Finally, take an arbitrary binary cubic $f(X,Y)$ with 
$\D_f\equiv \D\;(\mod\;k^{*6})$.  It is shown in 
[An, (3.4), (3.5), (3.8)] 
that the Jacobian of the curve $\C_f$ is the elliptic curve with 
affine model 
${Y^2=X^3-27\D_f}$, which is isomorphic to the curve 
$\E\:Y^2=X^3-27\D$ since 
$\D\equiv \D_f\;(\mod\;k^{*6})$.  	Thus, $\C_f$ is the curve of 
some 
$\p \in H^1(k,\E)$, and since $\C_f$ has a rational point in some 
degree 3 
extension of $k$, the restriction-corestriction argument shows t
hat $3\p=0$.  
Let ${\p_L=\res_{k\to L}(\p)\in H^1(L,\E)}$.  
	The earlier case with $\D$  a square shows that 
$\p_L=\Psi(\d)$ for some 
$\d\in H^1(L,\T)$.  This $\d$, which is not uniquely determined, 
need not 
lie in $\res_{k\to L}(H^1(k,\T))$.  But, let 
$\d'=2\cor_{L\to k}(\d)\in H^1(k,\T)$.  
	Then, using the compatibility of $\Psi$ with $\res_{k\to L}$ and 
$\cor_{L\to k}$, we have
$$
\align 
\res_{k\to L}(\Psi(\d'))  \ &= \  
\res_{k\to L}\big(2\cor_{L\to k}(\Psi_L(\d))\big) \ = \  
\res_{k\to L}(2\cor_{L\to k}(\p_L)) \\ 
&= \ \res_{k\to L}(2\dd[L\:\!k]\,\p) \ 
= \ 4\,\res_{k\to L}(\p) \ = \ \res_{k\to L}(\p). 
\endalign 
$$
	Because $2=[L\:\!k]$ is relatively prime to 3, $\res_{k\to L}$ is 
injective on the 3-torsion of $H^1(k,\E)$ and therefore 
$\Psi(\d')=\p$, showing 
that the curve $\C_f$ lies in $\im(\Psi)$.   \qed \enddemo

	Now, take any curve $\C_f$ as in (6.1) with $\D_f\not\in k^{*2}$, 
i.e., 
$f$ is not diagonalizable over $k$.  Set $\D=\D_f$ and take the 
elliptic 
curve $\E\: Y^2=X^3+c$ where $c=-27\D$.  
	As in \S 2 and \S 3, let 
$\T=\big\{\,\O,\;(0,\pm\sc)\,\big\} \subseteq \E(k_s)$, let $\E'$ be 
the elliptic curve $\E'\: Y^2=X^3+d$ where $d=-27c$ and let
$\T' =\big\{\,\O',\;(0,\pm\sqrt d)\,\big\} \subseteq \E'(k_s)$, and 
$\l'\:\E'\to \E$ the isogeny with kernel $\T'$ as in (3.6).  
	By Theorem 6.1, the curve $\E$ is the Jacobian of $\C_f$ and 
there is 
$\d\in H^1(k,\T)$ with $\Psi(\d)\in H^1(k,\E)$ corresponding 
to $\C_f$.  
Let $\partial\: \E(k)\to H^1(k,\T')$ be the connecting homomorphism 
arising from the short exact sequence~(3.2).  Since 
$\T'\otimes \T\cong \mu_3$ 
as $G_k$-modules, we have the cup product pairing  
$$
\cupp\: H^1(k,\T')\times H^1(k,\T)  \, \ra  \, H^2(k,\mu_3)  \, 
\cong  \, {}_3\Br(k),
$$ 
which yields a map   
$$
\b\:\E(k)\to {}_3\Br(k) \ \ \   \ \text{given by} \ \   \ 
P\mapsto \partial(P)\cupp \d.
\tag 6.6
$$  
The theorem of Ciperiani and Krashen [CK, Th.~2.6.5] says that 
$\Br(k(\C_f)/k)=\im(\b)$.

	As observed earlier in (2.4) and (3.7), we have 
$$
\T \, \cong  \, \zz_3(c) \, \cong \,  \mu_3(\D) \ \ \text{and} 
\ \ \T' \, \cong  \, \zz_3(\D) \, \cong  \, \mu_3(c) \ \ 
\text{as $G_k$-modules}.
$$
	In the diagonal case considered in earlier sections, when 
$\D\in k^{*2}$, 
the factors in the cup products  $H^1(k,\T')\cupp H^1(k,\T)
= H^1(k,\zz_3)\cupp H^1(k,\mu_3)$ give  
presentations of the corresponding algebras as cyclic algebras.  
However, when $\D\not\in k^{*2}$ and $-3\D\not\in k^{*2}$,  the cup 
products do not immediately yield cyclic algebra presentations, 
though 
we know by Wedderburn's theorem that these degree~3 algebras 
are cyclic.  
In this case, to realize the cup products as cyclic algebras one 
can apply 
the explicit algorithm given in the proof of [HKRT, Prop.~28] to 
restate 
a cup product in $H^1(k,\mu_3(c))\cupp H^1(k,\mu_3(\D))$  as a 
cup product in $H^1(k,\mu_3)\cupp H^1(k,\zz_3)$.
 	
	The cohomology groups $H^1(k,\T)$ and $H^1(k,\T')$ appearing here 
classify certain cubic field extensions of $k$, and we want to 
demonstrate 
how these field show up in the cup product algebras. 

\proclaim{Proposition 6.2} {\rm{(cf.\;[HKRT, Prop.~24])}} \ 
	Let $t\in k^*$.  Then the nontrivial cyclic subgroups of 
$H^1(k, \zz_3(t))$ classify $k$-isomorphism classes of separable 
field extensions $M$ 
of $k$ with $[M\:\!k]=3$ and $\disc(M)\equiv t\;(\mod\; k^{*2})$.
\endproclaim
\demo\nofrills{\smc {Proof.}\usualspace} \         
	(This is valid for any field $k$ with $\char(k)\ne2$.) \ This is 
standard when $t\in  k^{*2}$, so that $\zz_3(t)\cong \zz$ as 
$G_k$-modules.  In that case
for each nonzero $\chi\in H^1(k,\zz_3(t))\cong \Hom(G_k,\zz_3)$ 
the cyclic 
group $\langle\chi\rangle$ of order 3 corresponds to the fixed field 
$k_s^{\fker(\chi)}$, which is a cyclic Galois extension of $k$ of 
degree 3, 
so of trivial discriminant in $k^*/k^{*2}$.  All of the separable 
cubic field 
extensions of~$k$ with trivial discriminant arise in this way.  

	Now, suppose $t\not\in  k^{*2}$.  Let $K=k(\sqrt{t})$ and let $\s$ 
be the 
nonidentity $k$-automorphism of $K$.  The restriction map 
$H^1(k,\zz_3(t))\to 
H^1(K,\zz_3(t))= H^1(k,\zz_3)$ is injective as $\zz_3(t)$ is 
3-torsion and 
$[K\:\!k]=2$.  We identify $H^1(k,\zz_3(t))$ with its image 
in $H^1(K,\zz_3)$ 
which by Proposition 1.1 above consists of those 
$\chi\in H^1(K,\zz_3)$ with 
$\s(\chi)=\chi^{-1}$.  Let $N_{\chi}$ be the fixed field 
$k_s^{\fker(\chi)}$, 
a cyclic Galois extension of $K$ of degree 3.  Recall that 
$\s(\chi)$ is defined 
as follows:  Take any extension $\s'$ of $\s$ to $k_s$, so 
$\s'\in G_k$.  
For $\r\in G_K$, we then have $\s(\chi)(\r)=\chi(\s^{\prime -1}
\r\s')$.  Since $\s(\langle\chi\rangle)=\langle\chi\rangle$, it 
follows that 
$\ker(\chi)$ is a normal subgroup of $G_k$; hence $N_{\chi}$ is 
Galois over $k$.  
But the Galois group $\Gal(N_{\chi}/k)$ is nonabelian as 
$\s(\chi)\ne\chi$, so 
$\Gal(N_{\chi}/k)\cong S_3$, and $N_{\chi}$ contains three 
different but 
$k$-isomorphic subfields $M_i$ with $[M_i\:\!k]=3$.  Since 
$N_{\chi}$ is the 
normal closure of each $M_i$ over $k$, the discriminant of 
$M_i$ over $k$ is 
$tk^{*2}$ in $k^*/k^{*2}$.  The correspondence of the proposition 
is defined by 
mapping $\langle\chi\rangle$ to the isomorphism class of the 
$M_i$.  For the 
inverse of this map, consider a field extension $M$ of $k$ of 
degree 3 with 
$\disc(M)=tk^{*2}$.  Let $N=Mk(\sqrt t)$,  which is the normal 
closure of $M$ 
over $k$; so $N$ is Galois over $k$ with group $\Gal(N/k)\cong S_3$.  
Since $N$ 
is cyclic Galois over $k(\sqrt t)=K$, there is 
$\eta\in \Hom(G_K,\zz_3)$ with 
$G_N=\ker(\eta)$.  Because $N$ is Galois over $k$, 
$\s(\langle \eta\rangle)=\langle \eta\rangle$, but 
$\s(\eta)\ne \eta$  as 
$\Gal(N/k)$ is nonabelian; hence, $\s(\eta)=\eta^{-1}$.  The map 
$[M]\mapsto \langle \eta \rangle$ is an inverse of the map of 
this Proposition, so that  map is a bijection.    \qed\enddemo

\noindent{\bf Remark 6.3}.  \ Suppose $\mu_3\not\subseteq k$, and 
take any 
$u\in k^*\setminus k^{*2}$.  Then $H^1(k,\mu_3(u))$ injects into 
$H^1(k(\sqrt u),\mu_3)\cong k(\sqrt u)^*/k(\sqrt u)^{*3}$, but also 
$H^1(k,\mu_3(u))\cong H^1(k,\zz_3(-3u))$ which we have just seen 
classifies 
certain field extensions of $k$.  We can relate these two 
interpretations of 
$H^1(k,\mu_3(u))$ as follows:  Assume 
$u\not\equiv -3\;(\mod\;k^{*2})$.  Let 
$\r$ be the nonidentity $k$-automorphism of $k(\sqrt u)$.  Take 
any nonzero 
$\g\in H^1(k,\mu_3(u))$.  
	The injection $H^1(k,\mu_3(u))\to k(\sqrt u)/k(\sqrt u)^{*3}$ 
maps $\g$ to a 
cube class $sk(\sqrt u)^{*3}$ with $s\not\in k(\sqrt u)^{*3}$ 
such that 
$\r(s)\equiv s^{-1}\;(\mod\; k(\sqrt u)^{*3})$, say 
$\r(s)=s^{-1}b^3$.  Then 
$\r(b^3)=b^3$, so $\r(b)=b$ as $\mu_3\not\subseteq k(\sqrt u)$.  
Choose 
$r\in k_s$ with $r^3=s$, and let $F=k(\sqrt u)(r)$.  Since 
$\mu_3\not\subseteq k(\sqrt u)$, there is a unique extension 
of $\r$ to an 
automorphism $\r'$ of $F$, determined by $\r'(r)=r^{-1}b$.  
Since 
$\r^{\prime 2}(r)=rb^{-1}\r'(b)=r$, it follows that 
$\r^{\prime 2}=\id_F$.
 	Let $M$ be the fixed field~ $F^{\r'}$\!\!.  
Then $[M\:\!k]=3$, and one 
can check that $M$ has discriminant $-3uk^{*2}$ in 
$k^*/k^{*2}$ and that 
the isomorphism class of $M$ is the one associated to 
$\g\in H^1(k,\mu_3(-3u))$ 
in Proposition 6.2.

\medskip

\noindent{\bf Remark 6.4}.
	The degree 3 algebras $A$ in $\Br(k(\C_f)/k)$ with $f$ nondiagonal 
are 
realized as cup products of cohomology classes which are 
associated via 
Proposition 6.2 to certain cubic field extensions of $k$.  
It is interesting 
to see how those extension fields show up within $A$.  Take 
${\D\in k^*\setminus k^{*2}}$ with~$-3\D\not\in k^{*2}$.  Take any 
nonzero $\th\in H^1(k,\zz_3(\D))$ and $\g\in H^1(k,\mu_3(\D))$, 
and let 
$A$ be the degree 3 central simple $k$-algebra corresponding to 
$\th\cupp\g$ in ${}_3\Br(k)$.  Let  ${L=k(\sqrt \D)}$~and~let~ 
$\tau$~be~the~ 
non-identity $k$-automorphism of $L$.	Let $A_L=A\ok L$ and let 
${\th_L=\res_{k\to L}(\th)\in H^1(L,\zz_3)=\Hom(G_L,\zz_3)}$, and 
$\g_L=\res_{k\to L}(\g)\in H^1(L,\mu_3)\cong L^*/L^{*3}$.  Since 
$A_L=\th_L\cupp \g_L\in {}_3\Br(L)$, we can write 
${A_L=(N/L,\s,t)}$ where 
$N=k_s^{\text{\eightsl{ker}}(\th_L)}$, which is the cyclic 
Galois extension of $L$ of degree $3$ 
associated with $\th_L$, and $tL^{*3}\in L^*/L^{*3}$ corresponds 
to $\g_L$.  Because 
$\th_L\in \im\big(\res_{k\to L}H^1(k,\zz_3(\D))\big)$, 
by Proposition 1.1, $\tau(\th_L)=\th_L^{-1}$.  Consequently, 
$N$ is Galois but 
not abelian Galois over $k$, so ${\Gal(N/k)\cong S_3}$.  The 
automorphism 
$\tau$ of $L$ extends to an automorphism $\ol\tau$ of $N$ of order 2.  
(There are three different possibilities for $\ol\tau$; choose any 
one of 
them.) Then the fixed field $M=N^{\ol\tau}$ is a cubic extension of 
$k$ in the 
isomorphism class associated to $\th$ in Proposition 6.2.  
Note also that 
$\ol\tau\s\ol\tau^{\,-1}=\s^{-1}$ as $\Gal(N/k)\cong S_3$.  Now, 
$A_L= \ds{i=0}^{2} Nx^i$ where $xrx^{-1}=\s(r)$ for all $r\in N$ 
and $x^3=t$.  
	Because $\g_L\in \im\big(\res_{k\to L}(H^1(k,\mu_3(\D))\big)$,
 we have 
$\tau(t)=t^{-1}b^3$ for some $b\in L^*$.  Then $\tau(b^3)=b^3$ 
and hence 
$\tau(b)=b$ as $\mu_3\not\subseteq L$.   Extend $\ol\tau$ on $N$ to 
a map 
$\tt\:A_L\to A_L$ by sending $x\mapsto x^{-1}b$, i.e., 
$\tt(\sgm{i=0}^{2}r_ix^i)=\sgm{i=0}^{2}\ol\tau(r_i)(x^{-1}b)^i$.  
Clearly, $\tt$ is bijective.  
	Because $\tt$ respects the relations defining $A_L$ it is a ring 
isomorphism, which restricts to $\tau$ on $L$.  Also, as 
$\ol\tau^2=\id_N$ 
and $\tt^2(x)=xb^{-1}\tau(b)=x$, we have $\tt^2=\id$.  Thus, 
$\tt$ yields 
a semilinear action of $\Gal(L/k)$ on~$A_L$.
	Consequently, if we let $B=A_L^{\tt}$, the fixed ring of 
$A_L$ under the 
action of $\tt$, then $B$ is a $k$-algebra with $B\ok L\cong A_L$ 
($L$-algebra isomorphism) (cf.\;[J, p.\,56, Lemma 2.13.1]).  
Therefore, $B$~is a central simple $k$-algebra.   Since 
$\res\: {}_3\Br(k)\to {}_3\Br(L)$ is injective as $[L\:\!k]=2$, 
it follows 
that $B\cong A$.  Notice that $B$~contains  $N^{\tt}=N^{\ol\tau}=M$, 
the field 
corresponding to $\th$.  But, $\tt$ also maps $k(x)$ to itself 
and, as noted 
in Remark 6.3, the fixed field $k(x)^{\tt}$ is the degree 3 field 
extension 
of $k$ with discriminant $-3\D k^{*2}\in k^*/k^{*2}$ which 
corresponds to 
$\g$ viewed in $H^1(k,\zz_3(-3\D))$.  Thus, as~$A\cong B$, $A$ 
contains 
copies of the field $M$ corresponding to $\th$ and the field 
$F(x)^{\tt}$ 
corresponding to $\g$, in such a way that $M\ok L$ and 
$F(x)^{\tt}\ok L$ are 
the usual subfields in $A\ok L$ viewed as $\th_L\cupp\g_L$.

\medskip

	We now  give a few specific examples of relative Brauer 
groups of 
function fields of nondiagonal cubic forms over $\qq$. 

\noindent{\bf Example 6.5}.  \ Consider  $\C_f\: Z^3=f(X,Y)$ where 
$$
f(X,Y)\,=\,\, 13X^3+3\dd7X^2Y+3\dd19XY^2-5Y^3.
$$   The discriminant of $f$ is $88209=297^2 \in \qq^{*2}$, so $f$ is 
diagonalizable.
	Put $X=\frac13(U+2V)$ and $Y=\frac13(-U+V)$, which gives 
$\tf(U,V)=2U^3+11V^3$.   
Since $\C_f\cong \C_{\tf}$, from Example 5.13 we have 
$$
\Br(\qq(\C_f)/\qq)\,  =\, \Br(k(\C_f)/k) \, =\, \{\,0\,\},
$$    
where $k=\qq(\w)$.

\bigskip

\noindent{\bf Example 6.6}.  \ Consider $\C_f\: Z^3=f(X,Y)$ where 
$$
f(X,Y)\, =\, \, 4X^3+3\dd4XY^2+4Y^3.
$$ 
The discriminant $\D$ of $f$ is $320=2^6\dd5 \notin \qq^{*2}$, 
so $f$ is not 
diagonalizable over $\qq$ but is diagonalizable over 
$\qq(\sqrt 5\,)$. 
	The Jacobian of $\C$ is $ Y^2=X^3-27\D=X^3-27\dd2^6\dd5$.  
Thus, we can use the isomorphic elliptic curve 
$\E\:Y^2=X^3-27\dd5$, so $c = -27\dd5$,  and the isogenous 
curve ${\E'\:Y^2=X^3+5}$.  
By SAGE, it is 
known that $\E(\qq)$ has rank 1 with generator $(6,9)$ and that 
$\E'(\qq)$ has  generator $(-1,2)$.  
Note that $\a(6,9) =\ol{\,9+\sqrt {c\,}\,} = 
\ol{\,9+3\sqrt{-15}\,}$ in 
$\qq(\sqrt{-15})^*\big/\qq(\sqrt{-15})^{*3}$
and $\a'(1,2)= \ol{\,2+\sqrt 5\,} = \ol 1$ in~ 
$\qq(\sqrt 5)^*/\qq(\sqrt 5)^{*3}$, as  
$2+\sqrt 5 = \big((1+\sqrt 5)/2\big)^3$.
According to (5.7), we see that $\E(\qq)$ 
and $\E'(\qq)$ are 
torsion-free.   Thus, it follows from (5.4) that 
	$$
\a(\E(\qq))\, =\, \big\langle\, \ol{9+3\sqrt{-3\dd 5}\,}\,\big\rangle 
\subseteq \, \qq(\sqrt{-15})^*\big/\qq(\sqrt{-15})^{*3} \ \, \ 
\text{and} \  \ \, \a'(\E'(\qq))\, =\, \langle\, \ol{1}\,\rangle.
$$  
 	Because the map $\E(\qq)\to \Br(\qq(\C_f)/\qq)$ is surjective, 
the relative Brauer group $\Br(\qq(\C_f)/\qq)$ is completely 
determined by what 
this map does to the generator $(6,9)$.

	Let $k=\qq(\sqrt 5)$ and $L=\qq(\sqrt 5, \w)$.  
We have the following commutative diagram: 
$$
\CD  
\E(L)    @>>>  \Br(L(\C_f)/L) @= \Br(L(\C_{\tf})/L)  \\
 @AA{ }A              @AA{ }A         @AA{ }A    \\
  \E(k)    @>>>  \Br(k(\C_f)/k) @= \Br(k(\C_{\tf})/k)   \\
 @AA{ }A              @AA{ }A          \\
  \E(\qq)  @>>>  \Br(\qq(\C_f)/\qq)   \\
\endCD    
$$
	Here all the vertical maps are injective and the horizontal maps 
are surjective.  
We know $f$ is diagonalizable over $k$.  Specifically, put 
$X=\frac12(1+\sqrt 5)U+\frac12(1-\sqrt 5)V$ and $Y=-U-V$; then, 
$\tf(U,V)\, =\, \ta\,U^3+ \tb\,V^3$, where 
$$\textstyle
\ta\, =\,f\big(\frac12(1+\sqrt 5),-1\big) \,=\, 10(1+\sqrt 5)  
 \ \ \text{ and } \ \ 
 \ \tb\, =\,f\big(\frac12(1-\sqrt 5),-1\big)\,= \,  10(1-\sqrt 5).
$$   
We have  $\C_f\cong \C_{\tf}$ over $k$ and the  Jacobian is up 
to isomorphism  
$\E\:Y^2=X^3-27\dd\sqrt 5^2$. Thus, by Remark 4.2(ii) the image of 
$(6,9)\in \E(k)$ in 
$\Br(L(\C_{\tf})/L)$ is the Brauer class of the symbol algebra 
$$
A:=\, (t,u)_\w,
$$ 
where $t= -\tb/\ta = -(1 - \sqrt 5)\big/(1+ \sqrt 5)
=(3-\sqrt 5)\big/2$ and 
$u = 3(3 + \sqrt{-15})$. 
For~this~$t$~the~
nonidentity automorphism of $\qq(\sqrt 5)$
maps $t$ to $t^{-1}$.  Hence ${(t) \in 
H^1(\qq, \T) = H^1(\qq,\mu_3(5)) \subseteq \qq(\sqrt5)^*
\big/\qq(\sqrt 5)^{*3}}$. Likewise, the nonidentity automorphism
of $\qq(\sqrt{-15})$ maps $u$ to $3(3-\sqrt{-15}) = 6^3u^{-1}$,
so\break ${(u) \in H^1(\qq, \T') = H^1(\qq, \mu_3(-15)) \subseteq
\qq(\sqrt{-15})^*\big/\qq(\sqrt{-15})^{*3}}$.  
The algebra we seek is the cup product 
$(u)\cupp (t)$ in $H^2(\qq, \mu_3) \cong {}_3\Br(\qq)$, which 
maps to $A$ in~ 
${}_3\Br(L)$.  We will describe this algebra and its inverse by 
their local invariants.
Observe that each nonidentity $\sigma \in \Gal(L/\qq)$ maps exactly 
two of $t, u, \w$
to their inverses in $L^*/L^{*3}$ while fixing the third element; 
hence, 
$\sigma$ maps $[A]$ to $[A]$ in $\Br(L)$.  Therefore, the local 
invariants 
for $[A]$ are the same at each prime of $L$ over a given prime 
of $\qq$. This assures
that $A$ is defined over $\qq$, as expected.   
 
\def\L2{\widehat{L}_{2}}
\def\Q2{\widehat{\qq}_{2}}

      We now check that the symbol algebra $A$ above is nonsplit 
over $L$ by 
looking at a 2-adic completion.  For the 2-adic completion $\Q2$, 
note that 
$\Q2(\sqrt 5)$ and $\Q2(\sqrt {-3})$ are each the unramified 
quadratic extension of 
$\Q2$.  On the other hand, the 2-adic valuation on $\qq$ splits 
over $\qq(\sqrt{-15})$ 
so it  splits over $L$; let $\L2$ be the completion of $L$ with 
respect to either of 
it $2$-adic valuations.  Then $\L2$ is the unramified quadratic 
extension of $\Q2$.  
Let $v\:\L2^{\,*}\to \zz$ denote the corresponding discrete 
valuation.  Since the 
ramification of $A \otimes_L \L2$
in $H^1(\ol L_2,\zz_3)\cong 
\ol L_2^*/\ol L_2^{*3}$ is  the 
residue of 
$(-1)^{v(t)v(u)}t^{v(u)}u^{-v(t)}$, we will verify that 
$$
v(t)=0, \ \, \ol{t}\not\in\ol L_2^{*3}, \ \, \text{and} \ 
\,v(u)=1  \ 
\text{or} \ 2  \ \, \text{depending on the choice of } v.
\tag 6.7
$$  
Hence, $A \otimes_L \L2$ has nontrivial ramification,
so is nonsplit;, 
hence, $A$ is nonsplit. 

	To verify (6.7),  first observe that 
$t=-\frac{1-\sqrt 5}2\cdot\frac{2}{1+\sqrt 5}$ is a 2-adic unit.  
Since the 
minimal polynomial of $t$ over $\qq$ is $X^2-3X+1$, it follows 
that $\ol t$ is a 
root of 
$X^2+X+1$ over $\ol{L}_{2}\cong\Bbb F_4$.  Hence,  $\ol t\ne \ol 1$ 
in 
$\ol{L}_2^*$,
  so $\ol{t}\not\in\ol{L}_{2}^{*3}$.
In order to show that $v(u) = v(3+\sqrt{-15})=1$ or $2$, let 
$a= \sqrt{-15}$  over 
$\L2$.  Since $(a+1)(a-1) = -16$ and $a+1 = (a-1)+2$, the only 
possibilities are 
$v(a+1) = 1$ and $v(a-1) = 3$, or vice versa.  When $v(a-1) = 3$, 
we have 
$v(3 + \sqrt{-15}) = v(4+(a-1)) = 2$.  On the other hand, when 
$v(a+1) = 3$, 
we have $v(3 + \sqrt{-15}) = v(2+(a+1)) = 1$, as claimed in (6.7).

 	Now, notice that the only other prime where $A$ is potentially 
nonsplit is 
the unique 3-adic prime of $L$ since $t$ and $u$ are valuation 
units for all 
prime spots above $p\ne 2,3$ over $\qq$.  Since the sum of the 
local invariants 
is $0$ in $\qq/\zz$, 
  the algebra $B$ corresponding to $(u)\cupp (t)$ over $\qq$ has 
either 
local invariant $\frac13$ at 2, hence $\frac23$ at 3, and~0 
everywhere else; or 
$\frac23$ at 2, hence $\frac13$ at 3, and 0~everywhere else. The 
second possibility
corresponds to the inverse of the first in $\Br(\qq)$.  
Which possibility occurs depends on the choices of $\sqrt {-3}$, 
$\sqrt{5}$, and $\sqrt{-15}$.  In either case,
$$
\Br(\qq(\C_f)/\qq)  \, =\, \big\langle\,[(u)\cupp (t)]\,\big\rangle 
\,\cong\, \zz/3\zz,
$$
and the two algebras over $\qq$ we have specified by  local 
invariants
are the two nonidentity elements of this group.

\bigskip\bigskip         

\vfill\eject


\noindent {\bf \S7. \ The generalized Clifford algebra of $f$}          \medskip


	As mentioned in the Introduction and in the proof of 
Proposition 4.1, $\Br(k(\C_f)/k)$ consists of the specializations 
of the generalized Clifford algebra $\afk$ associated to 
a binary cubic form ~$f$ over ~
$k$; when $\D_f$ and $-3$ are squares in $k$, a presentation of 
the quotient ring of $\afk$ as a symbol algebra yields 
immediately a description of 
the algebras of degree 3 in $\Br(k(\C_f)/k)$ as symbol algebras.  
We now relate the results of the preceding sections to the 
structure of 
the Clifford algebra and its ring of central quotients $q(\afk)$.
We denote $q(\afk)$ by $\bfk$.  
We will give an explicit description of $\bfk$~as 
\roster
\item a symbol algebra if $f$ is diagonal and $-3\in k^{*2}$;  
\item a cyclic algebra if $f$ is diagonal and $-3\not\in k^{*2}$; 
\item a cup product in $H^1(Z,\Cal T)\cupp H^1(Z,\Cal T')$ if $f$ 
is not diagonalizable over $k$, 
\endroster
where $Z=Z(\bfk)\cong k(\E)$.  We will further show that 
$\bfk$ contains a maximal subfield isomorphic to $k(\E')$.  
We will also show how $\bfk$ is related to the curve $\C_f$ in the 
Brauer group $\Br(\E)$ of the elliptic curve $\E$ which is the 
Jacobian of $\C_f$.

	We work first with $f$ diagonal over a base field $L$, which can 
be any infinite field with ${\char(L)\ne 2,3}$.  Let 
$f(X,Y)=aX^3+bY^3\in L[X,Y]$, for some $a,b\in L^*$.  Let 
$\afl$ be the generalized Clifford algebra of $f$ over $L$.  Thus,
	$$
\afl=L\{x,y\}  \ \ \text{with defining relations} \;   
(px+qy)^3=f(p,q) 
\; \text{for all} \; p,q\in L^*.
$$
	It is known (see[H$_1$]) that $\afl$ is an Azumaya algebra over its center 
$Z=Z(\afl)$ which is an integral domain.  Let $q(Z)$ be the 
quotient field of $Z$, and let 
$\bfl=q(\afl) \cong\afl\otimes_{\!{}_Z} q(Z)$.  It is easy to see 
that the defining relations for $\afl$ are equivalent to the basic 
relations 
$$
x^3 \, = \, a,\  \ x^2y+xyx+yx^2 \, = \, 0,  
\ \ xy^2+yxy+y^2x\, =\, 0,  \ \ y^3\, =\, b.   \tag 7.1 
$$
	Let $\D_f= a^2b^2/4$, the discriminant of $f$, and let 
$c=-27a^2b^2/4=-27\D_f$.  Thus, for the projective curve 
$\C_f\: Z^3=f(X,Y)$, the Jacobian of $\C_f$ is the projective 
elliptic curve $\E\:Y^2=Z^3+c$. (As usual, we are describing 
projective elliptic curves in terms of   their affine models.)

	Let $F=L(\w)$, where $\w$ is a fixed primitive cube root of unity 
(possibly $F=L$).  For the Clifford algebra $\aff$ of $f$ over 
$F$, we have $\aff\cong \afl\otimes_{\!{}_L}F$.  We recall the 
description of $\aff$ given in~[H$_1$]:  Let 
$$
z\, =\, yx-\w xy \ \ \  \text{and}  \ \ \  \z\, =\, yx-\w^2xy.
$$
	The following identities follow easily from the basic identities 
$(7.1)$  (see [H$_1$, Lemma 1.5]): 
$$
xz\, =\, \w zx,  \ \ yz\, =\, \w^2zy,  \ \ x\z\, =\, \w^2\z x,
\  \ y\z\, =\, \w\z y,  \ \ z\z\, =\, \z z.   \tag 7.2
$$
One more useful fact is that 
$$
(yx^{-1})^3\, = \, b/a.  \tag 7.3
$$
	This is not transparent since $x$ and $y$ don't commute; here is a 
proof using the basic identities $(7.1)$:
	$$
\align (yx^{-1})^3 \, &=\, a^{-3}y(x^2y)(x^2yx^2)\, =\, 
-a^{-3}y(xyx+yx^2)x^2yx^2\, =\, -a^{-2}(yxy+y^2x)yx^2 \\  
&\, =\, a^{-2}xy^2yx^2\, =\,  b/a    \endalign 
$$
	To fix a choice of $\sc$, set $\sc=\frac32(\w-\w^2)ab
=\frac32(2\w+1)ab$.  Let 
	$$
\textstyle \rf\, =\, z\z \ \ \  \text{and} \ \ 
\sf\, =\, \frac12(z^3+\z^3).
$$
It follows from the identities in (7.2) that $z\z, \,z^3, \,\z^3$ 
commute with $x$ and $y$; hence $\rf$ and $\sf$ lie in the center 
$Z(\aff)$. 

\proclaim{Lemma 7.1 }
$z^3=\sf+\sc$ and  $\z^3=\sf-\sc$.  Hence, $\rf^3=\sf^2-c$. 
\endproclaim

\demo\nofrills{\smc {Proof.}\usualspace}
These essential formulas were proved in  [H$_1$, p.\;1274] using a 
matrix representation of $\aff$.   For the convenience of the reader,
 we now give a proof directly from the identities in (7.1) and (7.2).
Let $\alpha = yx$ and $\beta = xy$, so $z = \a -\omega \b$ and
$\z = \a - \omega^2 \b$.  
  Also, using (7.1), we have 
$$
\alpha \beta \, = \, -y(xyx + yx^2) \, = \, -(yxy+y^2x)x \, = 
\, \b \a.
$$
Hence,
$$
\align
\rf \, &= \, z\z \, = \, (\a - \omega \b)(\a- \omega^2\b)
\, = \, \a^2+ \a\b + \b^2 \\
&=\, (yxy)x +xy^2x + \b^2 \, = \, \b^2 - y^2x^2; 
\endalign
$$
likewise, $\rf =\a^2 -x^2y^2$. 
Hence, 
$$
\b^3\, =\, (\rf +y^2x^2)\b\, = \, \rf\b + y^2x^3y\,=\,
\rf\b +ab, \ \ \text{and likewise}  \ \ \a^3 \,=\, \rf\a + ab.
$$
Thus,
$$
\align
2ab\,&=\, (\a^3-\rf\a) + (\b^3-\rf\b)  \, = \,
(\a^3 + \b^3) -\rf(\a+\b) \\
&= \,
(\a^2-\a \b +\b^2 -\rf)(\a + \b) \, = \, -2\a\b(\a+\b).
\endalign
$$
This yields
$$
\align
z^3 - \z^3\,&=\, (\a-\omega\b)^3 - (\a-\omega^2\b)^3 \,
= \, (-3\w+3\w^2)\a^2\b +(3\w^2-3\w)\a\b^2 \\
&=\, (3\w^2 - 3\w)\a\b(\a+\b) \,=\, -(3\w^2 - 3\w)ab \,=\, 2\sqrt c.
\endalign
$$
Since $z^3+\z^3 = 2\sf$ by definition, we thus have 
$z^3 = \sf+\sqrt c$, \ $\z^3 = \sf-\sqrt c$, and 
$\rf^3 = z^3 \z^3 = \sf^2 -c$.\qed
\enddemo
 
	It was shown in [H$_1$, Th.~1.1]  that 
$\aff$ is an Azumaya algebra of rank 9 over its center~$Z(\aff)$, 
and that
$Z(\aff)=F[\rf,\sf]$,  and this ring is isomorphic to the coordinate 
ring over $F$  of an affine piece of $\E$.  Thus, for the quotient 
field we have $q(Z(\aff))=F(\rf,\sf)\cong F(\E)$.   

\proclaim{Proposition 7.2}  
	With notation as above, for $f = aX^3+bY^3$, the ring of quotients $\bff$
of the Clifford algebra $\aff$ is a symbol algebra,
$$
\bff \, = \,  F(\rf,\sf)\{x,z\}\, \cong\, 
 \(a,\, \sf+\sc; F(\rf,\sf)\)_{\w} \, \cong \, 
\(b/a,\, \sf+\sc; F(\rf,\sf)\)_{\w}\, ,
$$ 
and $\bfl$ is a cyclic algebra,
$$
\bfl\, =\,  \(L(\rf,\sf)(u)/L(\rf,\sf),\,\rho,\,a\) \, \cong\, 
 \(L(\rf,\sf)(u)/L(\rf,\sf),\,\rho,\,b/a\),
$$
where $u=z+\z$, which has minimal polynomial 
$X^3-\rf X-2\sf$ over $L(\rf,\sf)$.
\endproclaim 

\demo\nofrills{\smc {Proof.}\usualspace}  
	The formulas for $\bff$ except the last were given in [H$_1$], and 
follow easily  from the facts noted above.  Specifically, since 
$x^3=a$, $z^3=\sf+\sc$, and $xz=\w zx$, it follows that 
${F(\rf,\sf)\{x,z\}\cong \(a,\, \sf+\sc; F(\rf,\sf)\)_{\w}}$, 
which by 
dimension count is all of $\bff$.  	Likewise, by (7.3) and (7.2)
we have 
$(yx^{-1})^3=b/a$ and  $(yx^{-1})z=\w z(yx^{-1})$, so 
$\bff\cong \(b/a,\, \sf+\sc; F(\rf,\sf)\)_{\w}$.
 
	Let $K=F(\rf,\sf)(z)$, which is a maximal subfield of $\bff$ and 
is a 
cyclic Galois extension of $F(\rf,\sf)$ as $\w\in F(\rf,\sf)$
and  $z^3\in F(\rf,\sf)$ but 
$z\not\in F(\rf,\sf)$ as $z$ is not central in $\aff$.  Let 
$\rho$ be the 
$F(\rf,\sf)$-automorphism of $K$ with $\rho(z)=\w z$.  
	Then, $\rho(\z)=\rho(rz^{-1})=\w^2\z$.  Thus, for 
$u=z+\z\in K \setminus F(\rf,\sf)$, we have $K=F(\rf,\sf)(u)$ 
since $[K\:\!F(\rf,\sf)]=3$.  Also, we have 
$\rho(u)=\w z+\w^2\z=xux^{-1}=(yx^{-1})u(yx^{-1})^{-1}$.
	Hence, conjugation by $x$ or $yx^{-1}$ induces $\rho$ on $K$, 
showing that 
$$
\bff \, \cong\,  \(F(\rf,\sf)(u)/F(\rf,\sf),\,\rho,\,a\) \, \cong\, 
 \(F(\rf,\sf)(u)/F(\rf,\sf),\,\rho,\,b/a\).
$$  
By expanding out $u^3=(z+\z)^3$ using $z\z=\rf$ and $z^3+\z^3=2s$, 
we find that $u^3-\rf u-2\sf=0$; this determines 
the minimal polynomial 
of $u$  over $F(\rf,\sf)$.  This completes the proof if $L=F$.

	Suppose now that $L\ne F$.  We argue by descent from $F$ to $L$.  
Let $\s$ be the nonidentity $L$-isomorphism of $F=L(\w)$.  Let 
$\s$ also denote the canonical extension $\id_{\bfl}\otimes \s$ of~ 
$\s$ to ${\bff=\bfl\otimes_L F}$.  Since
${\aff=\afl\otimes_L F}$, we have  $\afl$ is the fixed ring 
$\aff^{\ \ \  \ \s}$ and ${\bfl=\bff^{\ \ \ \  \s}}$.  	Furthermore,  
$\s(\w)=\w^2$, 
$\s(x)=x$, and $\s(y)=y$, so  $\s(z)=\z$, $\s(\z)=z$, $\s(\rf)=\rf$, 
and $\s(\sf)=\sf$.  We thus obtain 
$Z(\afl)=(Z(\aff))^{\s}=F[\rf,\sf]^{\s}=L[\rf,\sf]$.  Since $\s(u)=u$,
 we have $\s(K) = K$ and $K^\s=F(\rf,\sf)(u)^\s=L(\rf,\sf)(u)$, 
with $[K^\s\!:\!L(\rf,\sf)] = [K\!:\!F(\rf,\sf)] =3$; so $K^\s$
 is a maximal subfield of $\bfl$.  
	Also, $\s\rho(u)= \w z+\w^2\z=\rho(u)$, so $\rho(K^\s) =  K^\s$ 
and $\rho|_{K^\s}$ is an 
automorphism (of order~3) of $K^\s$.  Hence, $K^\s$ is a cyclic 
Galois extension of $F(\rf,\sf)^\s=L(\rf,\sf)$, with 
$\rho|_{\!_{K^\s}}$ a  generator of the Galois group.  
	Conjugation by $x$ (or by $yx^{-1}$) induces $\rho$ on 
$K$ and hence on $K^\s$, and $x$, $yx^{-1} \in \bfl$  with $x^3$ 
and $(yx^{-1})^3$ central. This yields the desired descriptions of 
$\bfl$ as a cyclic algebra.  The minimal polynomial of $u$ over 
$F(\rf,\sf)$ has coefficients in $L(\rf,\sf)$, so it is also the 
minimal polynomial of $u$ over $L(\rf,\sf)$.  \qed\enddemo

	From the description of $\bfl$ and $\bff$ as cyclic and symbol 
algebras, one can see immediately that their specializations at the 
points of $\E(L)$ and $\E(F)$ coincide with $\Br(L(\C_f)/L)$ and 
$\Br(F(\C_f)/F)$ as given in Propositions 4.1 and 4.4 above.

	Now let  $d=-27c$ and let  $\E'$ be the projective elliptic curve 
$\E'\:Y^2=X^3+d$. It is clear from Cor.~4.6 and
the exact sequence 
$\E'(k) \rA{\l'} \E(k) \rA{\partial} H^1(k,\T')$  that the points in 
$\E(k)$ in the image of $\E'(k)$ under $\l'$ have 
trivial image in $\Br(k(\C_f)/k)$; so at these points 
the  specializations of $\bfk$ must be split algebras.     
	The next two  propositions show how to see this from the structure 
of $\bfk$.  We first show this in the diagonal case with ground 
field $L$ and $f(X,Y)=aX^3+bY^3$. We retain the notation used in 
Proposition 7.2 and its proof.  In particular we let 
$K=F(\rf,\sf)(u)$, a maximal subfield of $\bff$.

\proclaim{Proposition 7.3}  
	For the maximal subfield $K=F(\rf,\sf)(u)$ of $\bff$
as above, we have 
$K\cong F(\E')$.  If $\w \notin L$, so $F = L(\w) \supsetneqq L$, 
then $K^\s$ is a maximal subfield of $\bfl$, isomorphic to $L(\E')$.  
\endproclaim 

\demo\nofrills{\smc {Proof.}\usualspace}  
	Let $\rf'=6\sc/(z-\z)$ and $\sf'=9\sc\,(z+\z)/(z-\z)$, both lying in
 $K=F(\rf,\sf)(z)$.  Since $z^3=\sf+\sc$ and $\z=\rf z^{-1}$ with 
$\sf^2=\rf^3+c$, the proof of Lemma 3.1 
(with $\rf,\, \sf,\,  z,\, \z$ \, replacing $x,\, y,\, t, \, 
\widehat t$\, ) shows that 
$\, {\sf'}^2={\rf'}^3+d$, \, $\rf=({\rf'}^3+4d)/9{\rf'}^2$ and 
$\sf=({\sf'}^3-9d\sf')/27{\rf'}^3$.  Also,\break $z+\z=2\sf'/{3\rf'}$.  
Thus,  
$$
K\, =\, F(\rf,\sf)(z)\, =\,  F(\rf,\sf)(z+\z)\, = \, 
F(\rf,\sf,z+\z,\rf',\sf')\, =\, F(\rf',\sf').
$$
Since $K$ has transcendence degree 1 over $F$, 
$\, F(\rf',\sf')\cong F(\E')$.

	Now, assume that $L\ne F$, and let $\s$ be as in the proof of 
Proposition 7.2.  Since ${\s(z)=\z}$\break and $\s(\sc)=-\sc$ (as 
$\sc=\frac32(\w-\w^2)ab$), we have $\s(\rf')=\rf'$ and 
$\s(\sf')=\sf'$.  Hence,\break ${K^\s=F(\rf',\sf')^\s=L(\rf',\sf')}$, and 
$L(\rf',\sf')\cong L(\E')$ for the same reason that 
$F(\rf',\sf')\cong F(\E')$.  \qed \enddemo  


	Note  that the proof of Proposition 7.3 shows that the map on 
function fields $L(\E)\to L(\E')$ arising from the isogeny 
$\l'\:\E'\to \E$ given in (3.6) coincides with the inclusion 
$L(\rf,\sf)\to K^\s$.  
	Whenever a point $(r,s)\in \E(L)$ equals $\l'(r',s')$ for some 
$(r',s')\in \E'(L)$, then 
$$
{\l'}^{-1}(r,s)\, =\, \{(r',s'), \, (r',s')\oplus (0,\sqrt{d}),  \,  
(r',s')\oplus (0,-\sqrt{d})\} \, \subseteq \E'(L).
$$
The specialization of $K^\s$ at $(r,s)$  is the direct 
sum of the field extensions of $L$  determined by the points in 
${\l'}^{-1}(r,s)$, 
which here is $L\oplus L\oplus L$.
  Thus, the specialization of $\bfl$ at $(r,s)$ must split since it 
contains the specialization of $K^\s$.

	We now prove analogues  to Propositions 7.2 and 7.3 when the 
form is not diagonalizable over the ground field.  This will be 
done by descent from the quadratic extension over which the form 
is diagonalizable.  For this, let $k$ be the ground field
($\char(k)\ne 2,3$), and 
let 
$$
\fh(X,Y) \ = \  AX^3+3BX^2Y+3CXY^2+DY^3  \,\in k[X,Y]. \tag 7.4
$$  
As in (6.2), 
let $R = AC-B^2$, $2S = AD-BC$, $T = BD-C^2$, and $\D_{\fh}=S^2-RT$,
 which is assumed nonzero.
	We assume $\fh$ is not diagonalizable over $k$, so 
$\D_{\fh}\not\in k^{*2}$, which implies that $RT\ne 0$.  Let~$\d$~ 
be a fixed square root of $\D_{\fh}$ and let $L=k(\d)$.  Let 
$f(X,Y)=aX^3+bY^3$ be the diagonalization of $\fh$ over $L$ that 
was constructed in \S6. 

	To carry out the descent from $L$ to $k$, we must see how 
$\bfhk$ sits within $\bfl$.  For this, recall the change of 
variables used in \S6.  Matrix notation is convenient.  We write  
$f\M XY$ for $f(X,Y)$, likewise for $\fh$.
	As in \S6, set $U=RX+(S+\d)Y$ and $V=RX+(S-\d)Y$, i.e., 
$\M UV=P\M XY$ where 
$P=  \left(\smallmatrix R& \,S+\d \\ R & \,S-\d  
\endsmallmatrix\right)$. 
	Thus, $\M XY=Q\M UV$ where $Q=P^{-1}=\frac1{\,2R\d\,}
\left(\smallmatrix \d-S & \,\d+S \\ R  & \,-R  
\endsmallmatrix\right).$  
Then $f$ is obtained from~$\fh$ by the rule 
$f\M UV=\fh\(X\M UV,Y\M UV\!\)$, that is, $f=\fh\circ Q$.  
	Hence, in $L$,    
$$
a  \, =\,  f\M 10 \, =\,  \fh\(Q\M 10\!\) \, =\, 
 \fh \M{(\d-S)/2R\d}{1/2\d}   \ \ 
\text{and}  \ \ \ b\, =\,  f\M 01\, =\, 
\fh\M {(\d+S)/2R\d}{-1/ 2\d}.   \tag 7.5
$$  
	We have $\afl=L\{x,y\}$ where $f\M pq = (px+qy)^3 = 
\( \hmat xy\M pq\!\)^3$ 
for all $p,q\in L^*$.  
Let ${\afhk=k\{\xh,\yh\}}$, where 
$\fh\M{p'}{q'}= \( \hmat{\xh}{\yh}
\M{p'}{q'}\!\)^3$ for all $p', q'\in k$.  
	For any $p',q'\in k$, let $\M pq=P\M {p'}{q'}$, so 
$\M {p'}{q'}=Q\M pq$.  Then, we have   
$$
\fh\M{p'}{q'}\, =\, \fh\(Q\M{p}{q}\!\) \, =\,  f\M pq \, 
=\,  \( \hmat xy\M pq\!\)^3 \, =\,  
\(\hmat xy P\M {p'}{q'}\!\)^3.
$$
	If we set 
$\hmat {\wt x}{\wt y} =
\hmat xy P$, then 
$\fh\M{p'}{q'} = \hmat {\wt x}{\wt y}
\M{p'}{q'}$ for all $p',q'\in k$.  
	From this, we see that we can view $\afhk$ as a subalgebra of 
$\afl$ by setting $\hmat {\xh}{\yh}=
\hmat {\wt x}{\wt y} = 
\hmat xy P$.  
Thus, $\hmat xy=
\hmat {\xh}{\yh} Q$; that is,
$$ 
\textstyle x\,=\,\frac{\,\d-S\,}{2R\d}\,\xh\,+\frac1{\,2\d\,}\,\yh 
\ \  \ \ 
\text{and} \ \ \ \ y\,=\,\frac{\,\d+S\,}
{\,2R\d\,}\,\xh\,-\frac1{2\d}\,\yh.   
\tag 7.6
$$
	Then $\afhl=L\{\xh,\yh\}=\afl$.  Let $\tau$ be the $k$-automorphism 
of $L$ with 
$\tau(\d)= -\d$.  We extend 
$\tau$ to an automorphism of $\afhl=\afl$ by setting $\tau(\xh)=\xh$ 
and $\tau(\yh)=\yh$, and extend $\tau$ further to $\bfl=q(\afl)$.  
Thus, $\afhk=\afl^{\ \ \ \,\tau}$ and $\bfhk=\bfl^{\ \ \ \,\tau}$.

	Now, put $\ch=-27\D_{\fh}$ and $\dh=-27\ch$, and let $\Eh$ and 
$\Ehp$ be the projective elliptic curves $\Eh\:Y^2=X^3+\ch$ and 
$\Ehp\!\:Y^2=X^3+\dh$.  Let
$$
\textstyle \Th\, =\, \big\{\(0,\sqrt{\ch\,}\), 
\(0,-\sqrt{\ch\,}\),\,  
\O\big\} \subseteq \Eh(k_s) \ \ \  \text{and} \ \  \ 
\Thp \, =\,  \big\{ \(0,\sqrt{\dh\,}\), 
\(0,-\sqrt{\dh\,}\),\,  \O \big\}
\subseteq \Ehp(k_s).
$$ 
	Recall from \S6 that $\D_{\fh}=\e^6\D_f$, where 
$\e=\det(Q)=-1/(2R\d)\in L^*$; so, $\tau(\e)=-\e$.  
Finally, in $\aff$, set 
$$
z'\, =\, \e z, \ \ \z'\, =\, \e\z, \ \ 
\rfh\, =\, z'\z'\, =\, \e^2\rf, \ \ \text{and} \ \ \ 
\sfh\, ={z'}^3+\sqrt{\ch\,}\, =\, \e^3\sf.  
$$

\proclaim{Proposition 7.4}  With notation as above, 
for $\fh$ as in $($7.4$)$ we have the 
following$:$
\roster
\item "{\rm (a)}"  $Z(\afhk)=k[\rfh,\sfh]$ and \; 
$\sfh^2=\rfh^3+\ch$.  
Also, $k(\rfh,\sfh)\cong k(\Eh)$.
\item "{\rm (b)}"  Let $\Kh=k(\rfh,\sfh)(\uh)$, where 
$\uh=\w z'+\w^2\z'$, which has  minimal polynomial 
$X^3-\rfh X-2\sfh$ over $k(\rfh,\sfh)$.  Then $\Kh$ is a maximal 
subfield of $\bfk$ and $\Kh\cong k(\Ehp)$.
\item "{\rm (c)}"  $\bfhk$ is the cup product of the class of   
$b/a$ in $H^1(k(\rfh,\sfh),\Th)$ with the class of $\Kh$ in 
$H^1(k(\rfh,\sfh),\Th')$.
\endroster
\endproclaim   

\demo\nofrills{\smc {Proof.}\usualspace}  
	(a) (This was proved in [H$_1$, Th.~1.1$'$]. We give a more detailed 
argument here,  since we need it for the rest of the proof.)
Let $F=L(\w)$ where $\w^3=1$, $\w\ne 1$.  For convenience, we 
assume that $F\ne L$.  (When $F=L$, the Proposition still holds, 
and the arguments are a subset of what is given here.)  We use the 
same notation as from Proposition 7.2 and 7.3 for things associated 
to $\bff$, i.e., $z,\z,\rf,\sf,\s,\rho, K$.  
	Extend $\tau$ from $L$ to $F$ by setting $\tau(\w)=\w$, and extend 
$\tau$ likewise from $\afl$ (resp.~$\bfl$) to $\aff$ (resp.~$\bff$).  
Since $\tau(\s)=-\s$, formulas (7.5) and (7.6) show that $\tau(a)=b$, $\tau(b)=a$, $\tau(x)=y$, and $\tau(y)=x$ in $\afl$.  	Furthermore, since $c= -27\D_f=-27\e^6\D_{\fh}$, we have $\tau(\sc)=-\sc$. Also, $\tau(z)=\tau(yx-\w yx)=-\w\z$ and $\tau(\z)=-\w^2 z$, so $\tau(\rf)=\rf$, while $\tau(\sf)=-\sf$.  Hence, as $\tau(\e)=-\e$, we have $\tau(z')=\w\z' z$, $\tau(\z')=\w^2 z'$, $\tau(\rfh)=\rfh$, and $\tau(\sfh)=\sfh$.  Since $Z(\afl)=L[\rf,\sf]=L[\rfh,\sfh]$, we have $Z(\afhk)=Z(\afl)^{\tau}=k[\rfh,\sfh]$.  
	Also, $\sfh^2-\rfh^3=\e^6(\sf^2-\rf^3)=\e^6c=\ch$.  Furthermore, 
since $L(\rfh,\sfh)=L(\rf,\sf)$ is not algebraic over $k$ but 
$L(\rfh,\sfh)$ is algebraic over $k(\rfh,\sfh)$, it follows that 
$\rfh$ must be transcendental over $k$, which implies that 
$k(\rfh,\sfh)\cong k(\Eh)$.   

  	(b) Let $\kt=k(\rfh,\sfh)$, 
$\Lt=L(\rfh,\sfh)=L(\rf,\sf)=
\kt\(\sqrt{\botsmash{\D_{\widehat f}}
}\,\)$, and  
$\Ft=F(\rfh,\sfh)=F(\rf,\sf)=\Lt(\w)$.  Recall that $K$ is the 
field $\Ft(z)$, which is cyclic Galois over $\Ft$, and 
$K^{\s}=\Lt(u)$ where $u=z+\z$, and $K^{\s}$ is cyclic Galois 
over $\Lt$.  So, $K^{\s}=\Lt(z'+\z')$.  The $\Ft$-automorphism 
$\rho$ of $K$ with $\rho(z)=\w z$, $\rho(\z)=\w^2\z$ also 
satisfies $\rho(z')=\w z'$ and $\rho(\z')=\w^2\z'$, since 
$\rho(\e)=\e$.  
	Of the three conjugates $z'+\z'$, $\w z'+\w^2\z'$, 
$\w^2 z'+\w\z'$ over $\Lt$, $\tau$ fixes $\uh=\w z'+\w^2\z'$, 
while transposing the other two.  Since $K^\s=\Lt(u)=\Lt(\uh)$, 
it follows that for $K^{\s,\tau}$ \(meaning $(K^\s)^{\tau}$\) 
we have ${K^{\s,\tau} = \Lt(\uh)^{\tau} = \kt(\uh)=\Kh}$.  Hence, 
$\Kh$ is a maximal subfield of $\bfhk=(\bff)^{\s,\tau}$\!.  
	Set ${\rfhp=6\sqrt{\ch\,}\big/(\w z'-\w^2\z')}$ and 
${ \ \sfhp=9\sqrt{\ch\,}(\w z'+\w^2\z')\big/(\w z'-\w^2\z')}$ 
in $K$.  Since ${\tau(\sqrt{\ch\,})=\s(\sqrt{\ch\,})=-\sqrt{\ch\,}}$ 
and ${\tau(\w z')=\s(\w z')=\w^2\z'}$, we have 
$\rfhp,\sfhp \in K^{\s,\tau}=\Kh$.  
	Because $(\w z')^3=\sfh+\sqrt{\ch\,}$ and ${(\w z')(\w^2\z')=\rfh}$ 
with $\,\sfh^2=\rfh^3+\ch$, the argument for Proposition 7.3 shows 
that ${\, \sfhp}^2={\rfhp}^3+\dh$ and 
${\Kh=k(\sfhp,\rfhp) \cong k(\E')}$.  
Also, as $(\w z')(\w^2\z')=\rfh$ and 
$(\w z')^3+(\w\z')^3={z'}^3+{\z'}^3=2\,\sfh$, we have 
$\uh^3=(\w z'+\w^2\z')^3=\rfh\,\uh+2\,\sfh$, which yields the 
minimal polynomial of $\uh$ over $\kt$. 

	(c) We have $L\Kh=K^\s$, which we have seen is cyclic Galois over 
$\Lt$ with ${\Gal(K^{\s}/L)=\langle\rho |_{K^\s}\rangle}$ where 
$\rho$ is the $\Ft$-automorphism of $K$ with $\rho(z)=\w z$ and 
$\rho(\z)=\w^2\z$.  Note that\break 
${\rho\tau(z)=\rho(-\w\z)=-\z=\tau\rho^2(z)}$, so 
$\rho\tau=\tau\rho^2$.  So for the character 
$\chi\in H^1(\Lt,\zz_3)= \Hom(G_{\Lt},\zz_3)$ with kernel 
$G_{K^\s}$ and mapping $\rho$ to 1 (\mod\;3), we have 
$\(\tau(\chi)\)(\rho) = \chi(\tau^{-1}\rho\tau)=
\chi(\rho^{-1})= -\chi(\rho)$; thus, $\tau(\chi)=\chi^{-1}$.
	This implies by Proposition 1.1 that $\chi$ is the image of 
some $\p\in H^1(\kt,\Th')$ under the injective composition 
$H^1(\kt,\Th')\cong H^1(\kt,\zz_3(\D_{\fh}))\rA{\fres} 
H^1(\Lt,\zz_3)$.  Similarly, since ${\tau(b/a)=(b/a)^{-1}}$, the 
image of $b/a$ in $\Lt^*/\Lt^{*3}$ is the image of some 
$\psi\in H^1(\kt,\Th)$ under the injective map ${H^1(\kt,\Th') 
\cong H^1(\kt,\mu_3(\D_{\fh})) \rA{\fres} H^1(\Lt,\mu_3) 
\cong \Lt^*/\Lt^{*3}}$.
	Now, the map $\Br(\kt)\to \Br(\Lt)$ is injective on 
3-torsion subgroups as $[\Lt\:\! \kt]=2$.  In order to see that 
$\bfhk$ is the cup product of $\p$ with $\psi$, it suffices to 
check that $\bfhl$ is the cup product of $\res_{\kt\to\Lt}(\p)$ 
with $\res_{\kt\to\Lt}(\psi)$.  
	But, this is clear from the isomorphism 
$\bfhl=\bfl \cong (K^\s/L, \rho,b/a)$ of Proposition 7.2, as  
$\res_{\kt\to\Lt}(\p)$~is~$\chi$, whose kernel has fixed field 
$K^\s$, and $\res_{\kt\to\Lt}(\psi)=b/a \; (\mod\;L^{*3})$.    
\qed \enddemo

	We close with a remark on how $\bfk$ sits in $\Br(\E)$.  To 
simplify notation, let $f(X,Y)$ be an arbitrary nondegenerate 
binary cubic form over a field $k$, and let $\E$ be the Jacobian 
of the projective curve $\C_f\: Z^3=f(X,Y)$.  As we have seen, 
the ring of quotients $Q=\bfk$ of the Clifford algebra~$\afk$ has 
center which can be identified with the function field of the 
projective elliptic curve $\E$.  
	But more is known:  $Q$ is unramified at all the finite points 
of $\E$ because $\afk$ is an Azumaya algebra over the coordinate 
ring of the corresponding affine piece of $\E$; 
moreover, $Q$ is also 
unramified at the infinite point $\O$ of $\E$, and its 
specialization at $\O$ is split, as shown in [H$_1$, p.~1275].  
Therefore, $[Q]$ lies in $\Br(\E)$ the Brauer group of the 
elliptic curve $\E$.  This $\Br(\E)$ has equivalent 
characterizations~ as:
(1)~   the subgroup of $\Br(k(\E))$ consisting of the classes of 
algebras unramified at all points of $\E$ (not just rational 
points);  
(2) \  the group of equivalence classes of sheaves of Azumaya 
algebras over the curve $\E$; 
(3) \  $H^2_{\text{\'et}} (\E,\Bbb G_m)$.   
There is a more explicit description of $\Br(\E)$, as 
$\E$ is a projective elliptic curve:  Let $WC(\E)=H^1(k,\E)$, the 
Weil-Ch\^atelet group of $\E$, which classifies isomorphism 
classes of principal homogeneous spaces of $\E$ 
(cf.~[Si, pp.~290-291]).  
	There is a canonical short exact sequence 
$$
0 \ra \Br(k) \ra \Br(\E) \ra WC(\E) \ra 0,
$$ 
which is split by the specialization map at the rational point 
$\O$ which sends $\Br(\E)$ to $\Br(k)$.  (See~[Sk, p.~64] or 
[CK, (9)].  This is also deducible from the results in [L$_2$, \S 2].) 
	Thus, we have 
$$
\Br(\E) \, \cong  \, \Br(k)\oplus WC(\E).   \tag 7.7
$$
	Now, the curve $\C_f$ is a principal homogeneous space of $\E$, 
so it corresponds to some $\c_f\in WC(\E)$.  The class $[Q]$ in 
$\Br(\E)$ maps to $\c_f$ in $WC(\E)$, as was shown in 
[H$_3$, p.~518] 
though it was not stated that way.  Thus, $[Q]$ in $\Br(\E)$ 
corresponds to $(0,\c)$ in $\Br(k)\oplus WC(\E)$.  This 
illuminates the main result of [H$_3$], which says that 
$B$ is split if and only if $\C_f$ has a rational point.  
For, $C_f$ has a  rational point 
if and only if $C_f \cong \E$ over $k$ 
if and only if $\c_f=0$ in $WC(\E)$ 
if and only if $[Q]=0$ in $\Br(\E)$ by (7.7).

\bigskip\bigskip


\centerline{\bf REFERENCES}
\smallskip
\baselineskip=12pt
\def\hangbox to #1 #2{\vskip1pt\hangindent #1\noindent
\hbox to #1{#2}$\!\!$}

\def\myitem#1{\hangbox to 47pt {\hfill#1 \ \ }}

\frenchspacing
\bigskip

\myitem{[A]} A. A. Albert, On normal Kummer fields over a 
non-modular field, 
{\it Trans. Amer. Math. Soc.}, {\bf 36} (1934), 885-892.


\myitem{[An]} S. Y. An, S. Y. Kim, D. C. Marshall, S. H. Marshall, 
W. G. McCallum, and A. R. Perlis, Jacobians of genus one curves, 
{\it J.\;Number Theory}, {\bf 90} (2001), 304-315. 


\myitem{[Ca]} J. W. S. Cassels, {\it Lectures on Elliptic Curves}, 
Cambridge University Press, Cambridge, 1991.


\myitem{[CK]}   M. Ciperiani and D. Krashen, Relative Brauer groups 
of genus 1 curves, preprint, arXiv:math/0701614v2.


\myitem{[D]}   L. E. Dickson, {\it Algebraic Invariants}, J. Wiley, 
New York, 1914.


\myitem{[Gr]}  B. H. Gross, {\it Arithmetic on elliptic curves with 
complex multiplication}, Lecture Notes in Mathematics, {\bf 776}, 
Springer, Berlin, 1980

\myitem{[H$_1$]}  D.  Haile, On the Clifford algebra of a binary 
cubic form, 
{\it Amer.\;J.\;Math.}, {\bf 106} (1984), 1269-1280. 


\myitem{[H$_2$]}  D. Haile, On Clifford algebras, conjugate 
splittings, and 
function fields of curves, in {\it Israel Math.\;Conf.\;Proc}, 
volume 1: 
Ring Theory 1989, in honor of S.A. Amitsur (L.~Rowen, ed.) 
Weizman Science Press, 
Jerusalem, 356-361. 


\myitem{[H$_3$]}  D. Haile, When is the Clifford algebra of a 
binary cubic form split?, {\it J.\;Algebra}, {\bf 146} (1992), 
514-520. 


\myitem{[HKRT]} D.  Haile, M.-A. Knus, M. Rost, and J.-P. Tignol, 
Algebras of 
odd degree with involution, trace forms and dihedral extensions,  
{\it Israel J. Math.}, 
{\bf 96} (1996), 299-340. 


\myitem{[Hu]} D. Husem\"oller, {\it Elliptic Curves}, 2nd ed.,  
Springer, New York, 1987.


\myitem{[J]} N. Jacobson, {\it Finite Dimensional Division 
Algebras over Fields}, Springer, Berlin, 1996.
          

\myitem{[KMRT]} M.-A. Knus, A. S. Merkurjev, M. Rost, and 
J.-P. Tignol, {\it The Book of Involutions}, 
Amer.\;Math.\;Soc., Providence, RI, 1998.    


\myitem{[L$_1$]} S. Lichtenbaum, The period-index problem for 
elliptic curves,  
{\it Amer. J. Math.}, {\bf 90 } (1968), 1209-1223.
               

\myitem{[L$_2$]} S. Lichtenbaum, Duality theorems for curves 
over $p$-adic fields,  
{\it Invent. Math.}, {\bf 7} (1969), 120-136.


\myitem{[Mi]}  J. S. Milne, {\it Class Field Theory},   
http://www.jmilne.org/math/CourseNotes/CFT400.pdf 
 

\myitem{[SAGE]} {\it SAGE Mathematical Software}, Version 2.6, 
http://www.sagemath.org.


\myitem{[Se]} J.-P. Serre, {\it Local Fields}, Springer, 
New York, 1979.


\myitem{[Si]} J. Silverman, {\it The Arithmetic of Elliptic 
Curves}, Springer, New York, 1986.


\myitem{[Sk]} A. Skorobogatov, {\it Torsors and Rational Points}, 
Cambridge University Press, Cambridge, 2001.


\myitem{[V]} J. V\'elu, Isog\'enies entre courbes elliptiques, 
{\it C. R. Acad. Sci. Paris, S\'er. A-B}, {\bf 273} (1971), 
A238-A241.

\bigskip

\bigskip
$$
\openup-1.5\jot
\aligned
   &\text{Department of Mathematics}\\
&\text{Indiana University}\\
&\text{Bloomington, IN 47405}\\
&\\
&\text{e-mail: \tentt{haile\@indiana.edu }}
\endaligned\quad \ 
\aligned 
&\text{Department of Mathematics}\\ 
&\text{California State University}\\ 
&\text{San Bernardino, CA 92407}\\
&\\
&\text{e-mail: \tentt{ihan\@csusb.edu}}
\endaligned
\quad \ 
\aligned 
&\text{Department of Mathematics}\\ 
&\text{University of California, San Diego}\\ 
&\text{La Jolla, CA 92093-0112}\\
&\\ 
&\text{e-mail:
  \tentt{arwadsworth\@ucsd.edu}}
\endaligned
$$

\end